# RIEMANNIAN MANIFOLDS WITH MAXIMAL EIGENFUNCTION GROWTH

CHRISTOPHER D. SOGGE AND STEVE ZELDITCH

ABSTRACT. On any compact Riemannian manifold $(M, g)$ of dimension $n$, the $L^2$-normalized eigenfunctions $\{\phi_\lambda\}$ satisfy $||\phi_\lambda||_\infty \leq C\lambda^{\frac{n-1}{2}}$ where $-\Delta\phi_\lambda = \lambda^2\phi_\lambda$. The bound is sharp in the class of all $(M, g)$ since it is obtained by zonal spherical harmonics on the standard $n$-sphere $S^n$. But of course, it is not sharp for many Riemannian manifolds, e.g. flat tori $\mathbb{R}^n/\Gamma$. We say that $S^n$, but not $\mathbb{R}^n/\Gamma$, is a Riemannian manifold with maximal eigenfunction growth. The problem which motivates this paper is to determine the $(M, g)$ with maximal eigenfunction growth. Our main result is that such an $(M, g)$ must have a point $x$ where the set $\mathcal{L}_x$ of geodesic loops at $x$ has positive measure in $S_x^*M$. We show that if $(M, g)$ is real analytic, this puts topological restrictions on $M$, e.g. only $M = S^2$ (topologically) in dimension 2 can possess a real analytic metric of maximal eigenfunction growth. We further show that generic metrics on any $M$ fail to have maximal eigenfunction growth. In addition, we construct an example of $(M, g)$ for which $\mathcal{L}_x$ has positive measure for an open set of $x$ but which does not have maximal eigenfunction growth, thus disproving a naive converse to the main result.

## 1. INTRODUCTION

The problem which motivates this paper is to characterize compact Riemannian manifolds $(M, g)$ with maximal eigenfunction growth. Before stating our results, let us describe the problem precisely.

We first recall that the associated Laplace-Beltrami operator $\Delta = \Delta_g$, has eigenvalues $\{-\lambda_\nu^2\}$, where $0 \leq \lambda_0^2 \leq \lambda_1^2 \leq \lambda_2^2 \leq \ldots$ are counted with multiplicity. Let $\{\phi_\nu(x)\}$ be an associated orthonormal basis of $L^2$-normalized eigenfunctions. If $\lambda^2$ is in the spectrum of $-\Delta$, let $V_\lambda = \{\phi : \Delta\phi = -\lambda^2\phi\}$ denote the corresponding eigenspace. We then measure the eigenfunction growth rate in terms of

$$(1) \qquad L^\infty(\lambda, g) = \sup_{\substack{\phi \in V_\lambda \\ ||\phi||_{L^2}=1}} ||\phi||_{L^\infty}.$$

If $e_\lambda(f)$ is the projection of $f \in L^2$ onto $V_\lambda$, let $E_\lambda f = \sum_{\lambda_\nu \leq \lambda} e_{\lambda_\nu}(f)$ be the associated partial sum operators. Note that

$$(2) \qquad E_\lambda(x, y) = \sum_{\lambda_\nu \leq \lambda} \phi_\nu(x)\overline{\phi_\nu(y)}$$

Research partially supported by NSF grants.

is the kernel, where $\{\phi_\nu\}$ are as above. If

$$p(x, \xi) = \sqrt{\sum g^{jk}(x)\xi_j \xi_k} \tag{3}$$

is the principal symbol of $\sqrt{-\Delta}$, then by the local Weyl law [Av], [Le], [Ho 1],

$$E_\lambda(x, x) = \sum_{\lambda_\nu \leq \lambda} |\phi_\nu(x)|^2 = (2\pi)^{-n} \int_{p(x,\xi)\leq\lambda} d\xi + R(\lambda, x) \tag{4}$$

with uniform remainder bounds

$$|R(\lambda, x)| \leq C\lambda^{n-1}, \quad x \in M.$$

Since the integral in (4) is a continuous function of $\lambda$ and since the spectrum of the Laplacian is discrete, this immediately gives $\sum_{\lambda_\nu = \lambda} |\phi_\nu(x)|^2 \leq 2C\lambda^{n-1}$ (see Lemma 3.4), which in turn yields

$$L^\infty(\lambda, g) = 0(\lambda^{\frac{n-1}{2}}) \tag{5}$$

on any compact Riemannian manifold.

The bound (5) cannot be improved in the case of the standard round sphere, $(S^n, can)$, (by the zonal spherical harmonics) and on any rotationally invariant metric on $S^2$. It is however not attained by metrics of revolution on the 2-torus $T^2$, and the problem arises whether there is *any* metric $g$ on a surface of positive genus for which it is attained.

More generally, we pose:

- **Problem**: *Determine the $(M, g)$ for which $L^\infty(\lambda, g) = \Omega(\lambda^{\frac{n-1}{2}})$.*

Here we are using the notation that $\Omega(\lambda^{\frac{n-1}{2}})$ means $O(\lambda^{\frac{n-1}{2}})$ but not $o(\lambda^{\frac{n-1}{2}})$.

Our main result, Theorem 1.1, implies a necessary condition on a compact Riemannian manifolds $(M, g)$ with maximal eigenfunction growth: there must exist a point $x \in M$ for which the set

$$\mathcal{L}_x = \{\xi \in S_x^* M : \exists T : \exp_x T\xi = x\} \tag{6}$$

of directions of geodesic loops at $x$ has positive surface measure. Here, exp is the exponential map, and the measure $|\Omega|$ of a set $\Omega$ is the one induced by the metric $g_x$ on $T_x^* M$. For instance, the poles $x_N, x_S$ of a surface of revolution $(S^2, g)$ satisfy $|\mathcal{L}_x| = 2\pi$.

THEOREM **1.1**. *Suppose that $|\mathcal{L}_x| = 0$. Then given $\varepsilon > 0$ there exists a neighborhood $\mathcal{N} = \mathcal{N}(\varepsilon)$ of $x$, and a positive number $\Lambda = \Lambda(\varepsilon)$, so that*

$$\sup_{\phi \in V_\lambda} \frac{\|\phi\|_{L^\infty(\mathcal{N})}}{\|\phi\|_{L^2(M)}} \leq \varepsilon \lambda^{(n-1)/2}, \quad \lambda \in spec\ \sqrt{-\Delta} \geq \Lambda. \tag{7}$$

*If one has $|\mathcal{L}_x| = 0$ for every $x \in M$ then given $\varepsilon > 0$ there exists $\Lambda(\varepsilon, p)$, so that*

$$\sup_{\phi \in V_\lambda} \frac{\|\phi\|_{L^p(M)}}{\|\phi\|_{L^2(M)}} \leq \varepsilon \lambda^{\delta(p)}, \quad \lambda \in spec\ \sqrt{-\Delta} \geq \Lambda(\varepsilon, p), \quad p > \frac{2(n+1)}{n-1} \tag{8}$$

*where*

$$\delta(p) = \begin{cases} n(\frac{1}{2} - \frac{1}{p}) - \frac{1}{2}, & \frac{2(n+1)}{n-1} \leq p \leq \infty \\ \frac{n-1}{2}(\frac{1}{2} - \frac{1}{p}), & 2 \leq p \leq \frac{2(n+1)}{n-1}. \end{cases} \tag{9}$$

The $L^p$-bonds improve on some estimates in [So2], which say that

$$\sup_{\phi \in V_\lambda} \frac{\|\phi\|_p}{\|\phi\|_2} = O(\lambda^{\delta(p)}), \quad 2 \leq p \leq \infty. \tag{10}$$

The final result is false in the endpoint case $p = \frac{2(n+1)}{n-1}$, as will be explained in §8.

Recently, Taylor [Ta] proved bounds like (7) under more restrictive hypotheses involving the caustics of the wave group, and in [Ta] he also discusses applications of these sort of bounds to problems involving the convergence of eigenfunction expansions. Also, Sarnak [Sa] has recently obtained limit formulae for the $L^4$-norm of certain eigenfunctions on arithmetic hyperbolic quotients.

Our eigenfunction bounds are based on the following estimates for the local Weyl remainder.

THEOREM **1.2.** *As in (4) let $R(\lambda, x)$ denote the remainder for the local Weyl law at $x$. Then*

$$R(\lambda, x) = o(\lambda^{n-1}) \text{ if } |\mathcal{L}_x| = 0. \tag{11}$$

*Additionally, if $|\mathcal{L}_x| = 0$ then, given $\varepsilon > 0$, there is a neighborhood $\mathcal{N}$ of $x$ and a $\Lambda = <\infty$, both depending on $\varepsilon$ so that*

$$|R(\lambda, y)| \leq \varepsilon \lambda^{n-1}, \ y \in \mathcal{N}, \ \lambda \geq \Lambda. \tag{12}$$

As we were informed by V. Guillemin, A. Laptev and D. Robert after completing this paper, these estimates overlap previous results of Y. Safarov in the article [Saf] (see also [SV], §1.8 and [Iv2]). Because our proof is somewhat different from that of Safarov and our exposition contains complete details, we have decided to retain our original proof.

Theorem (1.2) is a kind of local analogue of the classical Duistermaat-Guillemin [DG] (see also Ivriĭ[Iv1]) result that $N(\lambda) = (2\pi)^{-n} \text{Vol}(M) \lambda^n + o(\lambda^{n-1})$ if the set of initial directions $(x, \xi)$ of closed geodesics has measure zero in $S^*M$. Here, $N(\lambda)$ denotes the number of eigenvalues $\leq \lambda$ (counted with multiplicity). Their assumption is the natural one for this theorem since only closed geodesics contribute singularities to the trace of the wave group $U(t)$, $t \neq 0$; however, for pointwise bounds at $x$ one needs to assume zero measure of closed loops through $x$.

We should point out that (11) by itself is not strong enough to give the above sup-norm estimates for eigenfunctions. They also require the more delicate second statement, which says that $x \to \limsup_{\lambda \to \infty} \lambda^{n-1} R(\lambda, x)$ is continuous at points where $|\mathcal{L}_x| = 0$. We now give some concrete applications and examples of these general results.

Our first application is to the real analytic setting, the result has strong implications on the topology of $M$.

THEOREM **1.3.** *Suppose that $(M, g)$ is real analytic and that $L^\infty(\lambda, g) = \Omega(\lambda^{(n-1)/2})$. Then $(M, g)$ is a $Y_\ell^m$-manifold, i.e. a pointed Riemannian manifold $(M, m, g)$ such that all geodesics issuing from the point $m$ return to $m$ at time $\ell$. In particular, if $\dim M = 2$, then $M$ is topologically a 2-sphere $S^2$.*

For the definition and properties of $Y_\ell^m$-manifolds, we refer to [Besse] (Chapter 7). By a theorem due to Bérard-Bergery (see [BB, Besse], Theorem 7. 37), $Y_\ell^m$ manifolds $M$ satisfy $\pi_1(M)$ is finite and $H^*(M, \mathbb{Q})$ is a truncated polynomial ring in one generator. This of course implies $M = S^2$ (topologically) when $n = 2$. We remark that the loops are not assumed to close up smoothly. An interesting example to keep in mind here is the tri-axial ellipsoid $E_{a_1,a_2,a_3} = \{(x_1, x_2, x_3) \in \mathbb{R}^2 : \frac{x_1^2}{a_1^2} + \frac{x_1^2}{a_1^2} + \frac{x_1^2}{a_1^2} = 1\}$, with $a_1 < a_2 < a_3$. This is a metric on $S^2$ with four umbilic points. At each umbilic point, all the geodesics leaving the point return to it at the same time. Only two such geodesics close up smoothly, namely the middle closed geodesic $\{x_2 = 0\}$, traversed in the two possible directions. For more on the geometry of the ellipsoid, we refer to [A, K, CVV].

Our second application is to generic metrics:

THEOREM 1.4. $L^\infty(\lambda, g) = o(\lambda^{(n-1)/2})$ *for a generic Riemannian metric on any manifold.*

The proof is just to show that the condition $|\mathcal{L}_x| = 0$, $\forall x \in M$ holds for a residual set of metrics with respect to the Whitney $C^\infty$ topology. It appears that this is a new geometric result.

To our knowledge, the most general class of metrics for which an improvement of the general error bound $O(\lambda^{\frac{n-1}{2}})$ in the local Weyl law has previously been proved is that of manifolds without conjugate points satisfying exponential bounds on the geodesic flow [Be]. In that case, the error term has been improved to $O(\frac{\lambda^{(n-1)/2}}{\log \lambda})$. The result is not stated explicitly in [Be], but follows from its estimates on the remainder term in the Weyl law. Estimates of the error term in the Weyl law also appear in [V], but it appears that the methods of that paper require integration in $x$ and hence do not imply $L^\infty$ bounds; in particular, their assumptions are on closed geodesics rather than loops.

It is natural to ask whether the converse of these results is also true. Regarding the first remainder estimate in Theorem 1.2, we shall prove the conditional converse:

THEOREM 1.5. *If the loopsets of $(M, g)$ are clean, or even almost-clean, then the following implication holds:*
$$R(\lambda, x) = o(\lambda^{n-1}) \implies |\mathcal{L}_x| = 0.$$

Clean and almost-clean loopsets are generalizations of the same notions for fixed point sets and will be defined in §7. Almost-clean $(M, g)$ include manifolds containing 'round annuli', i.e. annuli isometric to equatorial annuli of the round sphere. One can obtain such manifolds by deforming the round metric on $S^n$ in small balls, say by adding handles. The geodesic flow in such cases is not clean in the sense of Duistermaat-Guillemin [DG] but is 'usually' clean outside a hypersurface (the boundaries of the balls). The proof of Theorem 1.5 is similar to the proof of the converse to the Duistermaat-Guillemin theorem in [Z] and the assumption of almost-cleanness is taken from there.

However, the naive converse to our sup-norm result is simply false.

THEOREM 1.6. *There exist $C^\infty$-Riemannian tori of revolution $(T^2, g)$ such that:*
$$\exists x : |\mathcal{L}_x| > 0, \ R(\lambda, x) = \Omega(\lambda^{\frac{n-1}{2}}), \ but \ \frac{||\phi_\lambda||_\infty}{||\phi_\lambda||_2} = o(\lambda^{\frac{n-1}{2}}).$$

Thus, remainder blow-up is not sufficient for maximal eigenfunction growth. We note that it is very likely that the ellipsoid is also an example of such behaviour: we have $R(\lambda, x) = \Omega(\lambda^{\frac{n-1}{2}})$ at each umbilic point, but it is conjectured that $\frac{||\phi_\lambda||_\infty}{||\phi_\lambda||_2} = o(\lambda^{\frac{n-1}{2}})$ (J. Toth). In the final section, we shall discuss some interesting open problems on maximal eigenfunction growth.

This paper is organized as follows. In §2, we review some basic results on the wave group and geodesic flow, and on several Tauberian theorems. This background is sufficient for the proofs of Theorems 1.1 and 1.2 in §3. The proof of Theorem 1.2 is in part a variant of Ivriǐ's [Iv1] (see also [Ho IV]) proof of the Duistermaat-Guillemin theorem [DG] that the remainder term in Weyl's law is $o(\lambda^{\frac{n-1}{2}})$ if the set of closed geodesics has measure zero). In particular, we shall need to analyze the small-time behavior of the restriction to the diagonal in $M$ of the kernel of $CU(t)B$, if $C$ and $B$ are zero-order pseudodiferential operators. The case where $C = Id$ is the identity operator was analyzed by Ivriǐ[Iv1] [Ho IV]; however, since we cannot take traces in our proof, we need to study this slightly more technical case. Fortunately, we are able to reduce the analysis to the special case discussed in §2. Additionally, Theorem 1.2 uses some continuity properties of $|\mathcal{L}_x|$. The proof of Theorem 1.1 for the $L^\infty$ norms is then a simple matter, and the other cases follow by interpolating the $L^p$ bounds in [So1]. We then apply the results to some concrete classes of metrics. As preparation, we discuss the geometry of loops in §4. In §5 and §6, we give the applications of Theorem 1.1 to real analytic manifolds and generic Riemannian manifolds, respectively. In §7, we study converses to Theorems 1.1 and 1.2, and construct an example where $R(\lambda, x) = \Omega(\lambda^{\frac{n-1}{2}})$ but $L^\infty(\lambda, g) = o(\lambda^{\frac{n-1}{2}})$. Finally, in §8, we consider related open problems.

Regarding notation, $C$ will denote a constant which is not necessarily the same at each occurrence. Also, integrals like the one in (4) should be interpreted as $C^\infty$ densities which are homogeneous in $\lambda$ and defined so that the scalar product with $\psi \in C^\infty(M)$ is the integral of the product by $\psi$ with respect to the symplectic measure $dxd\xi$ on $T^*M$.

It is a pleasure to thank R. Hardt, W. Minicozzi, B. Shiffman for very helpful conversations on analytic sets, and for giving us copies of the unpublished articles [L, S]. We are particularly grateful to V. Guillemin, A. Laptev and D. Robert for the reference to Safarov's work [Saf]. We also thank J. Toth for bringing the ellipsoid example to our attention.

## 2. Wave group and geodesic flow

In this section, we collect the relevant background on the wave group and geodesic flow for the proofs of Theorems 1.2 and 1.1. Essentially, they are based on the short-time asymptotics of the wave equation, long time wave front relations and Tauberian theorems. We shall need substantially more geometry of loopsets to prove Theorems 1.5 and 1.6, but postpone the discussion until §4

Throughout this paper, $(M, g)$ denotes a compact Riemannian manifold of dimension $n \geq 2$, $T^*M$ denotes its cotangent bundle and $S^*M$ denotes its unit sphere bundle with respect to $g$. We denote by $\exp tH_p$ the geodesic flow of $g$, defined as the flow of the Hamiltonian vector field $H_p$ of $p(x, \xi) = \sqrt{\sum g^{jk}(x)\xi_j\xi_k}$, the principal symbol for $\sqrt{-\Delta}$. By definition, $\exp tH_p$ is homogeneous, i.e. commutes with the natural $\mathbb{R}^+$-action, $r \times (x, \xi) = (x, r\xi)$ on $T^*M\setminus 0$. We also define the exponential map at $x$ by $\exp_x \xi = \pi \circ \exp tH_p(x, \xi)$. These definitions are standard in microlocal analysis but differ from the usual geometer's definitions, which takes $p^2 = \sum_{i,j} g^{ij}(x)\xi_i\xi_j$ as the Hamiltonian generating the geodesic flow. The geometer's geodesic flow is not homogeneous.

2.1. **Wave group and spectral projections.** Let $dE_\lambda$ be the spectral measure associated with the partial summation operators $E_\lambda$ for $\sqrt{-\Delta}$. Then we shall prove our estimates for the local Weyl law remainder term $R(\lambda, x)$ by studying the Fourier transform

$$U(t) = \int e^{-it\lambda} dE_\lambda = e^{-it\sqrt{-\Delta}}.$$

This of course is the wave group, i.e. the unitary group generated by $\sqrt{-\Delta}$. It solves the Cauchy problem

$$(\frac{1}{i}\frac{\partial}{\partial t} + \sqrt{-\Delta})U(t) = 0, \quad U(0) = Id.$$

It is well known that the kernel $U(t, x, y)$ is a Fourier integral operator in the class $I^{-1/4}(\mathbb{R} \times M \times M, \mathcal{C})$ where $\mathcal{C}$ is the Lagrangean

$$\{(t, x, y; \tau, \xi, \eta) : \tau + p(x, \xi) = 0 \text{ and } (x, \xi) = \exp(tH_p(y, \eta))\}.$$

The last condition means that the orbit of $H_p$ which is at $(y, \eta)$ at time 0 is at $(x, \xi)$ at time $t$. In particular, $x$ and $y$ are joined by a geodesic of length $|t|$.

We shall mainly be concerned with the restriction of the wave kernel to the diagonal in $M \times M$. The restriction is defined by

$$U(t, x, x) = \Delta^* U(t, x, y), \quad \text{where } \Delta : \mathbb{R} \times M \to \mathbb{R} \times M \times M \text{ is } \Delta(t, x) = (t, x, x).$$

Although the notation conflicts with that for the Laplacian, it is standard and should cause no confusion. By general wave front set considerations, one always has

$$WF(U(t, x, x)) \subset \mathcal{C}_\Delta,$$

where

(13) $$C_\Delta = \{(t, \tau, x, \xi - \eta) : \tau = -|\xi|, \ \exp tH_p(x, \eta) = (x, \xi)\}.$$

Indeed, the canonical relation underlying $\Delta^*$ is given by ([DG], (1.20)),

$$WF'(\Delta) = \{(t, \tau, x, (\xi + \eta)); (t, \tau, x, \xi, x, \eta))\} \subset T^*((\mathbb{R} \times M) \times (\mathbb{R} \times M \times M)).$$

The intersection defining $C_\Delta$ is essentially that between $\mathcal{C}$ and the second component

$$\mathcal{D} := \{(t, \tau, x, \xi, x, \eta)\} \subset T^*(\mathbb{R} \times M \times M)$$

of $WF'(\Delta)$. If $\mathcal{C} \cap \mathcal{D}$ is a clean intersection (cf. [DG], [Ho IV], §??), $C_\Delta$ is a Lagrangian submanifold of $T^*(\mathbb{R} \times M)\backslash 0$, and one has the much stronger statement that

$$U(t, x, x) \in I^0(\mathbb{R} \times M, \mathcal{C}_\Delta),$$

where $I^0(\mathbb{R} \times M, \mathcal{C}_\Delta)$ is the class of Fourier integral operators of order zero associated to the canonical relation $C_\Delta$. However, we shall not need this hypothesis in the proofs of Theorems 1.2-1.1.

We note $(x, \eta)$ is the initial (co-)tangent vector and $(x, \xi)$ is the final (co-)tangent vector of a geodesic loop of length $t$ at $x$; loops need not be smoothly closed. If there are no geodesic loops of length $|t_0|$ through $x_0 \in M$, it follows that $(t, x) \to U(t, x, x)$ is free of singularities at $(t_0, x_0)$. Since the length of a geodesic loop is at least the injectivity radius $\delta$ of $(M, g)$, it follows that $U(t, x, x)$ is smooth for $|t| < \delta$.

2.2. **Small-time asymptotics for microlocal wave operators.** In the proof of Theorem 1.2, we shall need to calculate asymptotics for the small-time behavior of the diagonal restriction $[U(t)Q](x,x)$ $\{(x,x) : x \in M\}$ of the one-parameter family of operators

$$U(t)Q : C^\infty(M) \to C^\infty(M),$$

where $Q(t, x, D_x)$ is a polyhomogenous zero-order pseudodifferential operator. In the case where $Q$ is independent of $t$, this calculation was carried out by Hörmander [Ho IV], following Ivrii[Iv1], and, as we shall see, the slightly more general case mentioned above follows in a straightforward way from this special case.

To set up the notation, write

(14) $$Q(t, x, D_x) = B(x, D_x) + tR(t, x, D_x),$$

where of course $B(x, D_x) = Q(0, x, D_x)$. We then shall let $b$ and $b^s$ denote the principal and subprincipal symbols of $B$, while $r$ shall denote the principal symbol of $R(0, x, D_x)$. We recall that $p = \sqrt{\sum g^{jk}(x)\xi_j \xi_k}$ denotes the principal symbol of $\sqrt{-\Delta}$, and that $p^s = 0$ for this operator.

The calculation we require then is the following

**Proposition 2.1.** *For any compact Riemannian manifold $(M,g)$, the restriction $K$ of the kernel of $U(t)Q$ to $\mathbb{R} \times \{(x,x) : x \in M\}$ is conormal with respect to $\{0\} \times \{(x,x) : x \in M\}$ in a neighborhood of this submanifold. Moreover, there is a $\delta > 0$ so that when $|t| < \delta$*

(15) $$K(t,x) = \int_{-\infty}^{\infty} \frac{\partial A(x, \lambda)}{\partial \lambda} e^{-i\lambda t} \, d\lambda,$$

*where $A \in S^n$, $A(x, 0) = 0$ and*

(16) $$A(x, \lambda) - (2\pi)^{-n} \int_{p(x,\xi)<\lambda} (b + b^s) d\xi + (2\pi)^{-n} \frac{\partial}{\partial \lambda} \int_{p(x,\xi)<\lambda} \tfrac{i}{2}\{b, p\} \, d\xi$$
$$+ (2\pi)^{-n} \frac{\partial}{\partial \lambda} \int_{p(x,\xi)<\lambda} ir(0, x, \xi) \, d\xi \in S^{n-2}.$$

If $Q = B$ and thus is independent of $t$ then (16) is just Proposition 29.1.2 in [Ho IV]. Since the kernel $K(t, y)$ in this case is conormal, it is easy to see that this special case also yields the time-dependent case. The constant $\delta$ is the injectivity radius of $(M, g)$

Indeed if

(17) $$Q(t, x, D_x) = tR(t, x, D_x) \sim t \sum_{j=0}^{\infty} \frac{t^j}{j!} \frac{\partial^j R}{\partial t^j}(0, x, D_x)$$

then we can compute the contributions of the various terms of the expansion to (16). We first note that, by the aforementioned time-independent case, the $j = 0$ term will contribute a term of the form

$$K_0(t,x) = \int t \frac{\partial A_0}{\partial \lambda}(x, \lambda) e^{-i\lambda t} d\lambda = \int -i \frac{\partial^2 A_0}{\partial \lambda^2} e^{-i\lambda t} d\lambda$$

with $A_0 \in S^n$ satisfying

$$A_0(x,\lambda) - (2\pi)^{-n}\int_{p<\lambda}(r+r^s)d\xi + (2\pi)^{-n}\frac{\partial}{\partial\lambda}\int_{p<\lambda}\tfrac{i}{2}\{r,p\}d\xi \in S^{n-2},$$

if $r$ and $r^s$ are the principal and subprincipal symbols of $R(0,x,D_x)$.

Note that

$$-i\frac{\partial A_0}{\partial\lambda}(x,\lambda) + (2\pi)^{-n}\frac{\partial}{\partial\lambda}\int_{p<\lambda}ird\xi \in S^{n-2},$$

and so $K_0$ exactly contributes the last term in the right side of (16). Since this argument also shows that the remaining piece

$$\sim t\sum_{j>0}\frac{t^j}{j!}\frac{\partial^j R}{\partial t^j}(0,x,D_x)$$

only contributes a term $\partial\tilde{A}/\partial\lambda$ to (15) with $\tilde{A} \in S^{n-2}$, the proof is complete. □

To prove Theorem 1.2 we shall have to write

$$V(t) = CU(t)B$$

as an expansion like (15)-(16) if $C(x,D)$ and $B(x,D)$ are zero-order polyhomogeneous pseudodifferential operators on $M$. In Ivriĭ's [Iv1] proof of the Duistermaat-Guillemin [DG] theorem (see also [Ho IV]) it sufficed just to compute the asymptotics for $UB$ since the ultimate calculation involved a trace over $M$, rather than a pointwise estimate like the one we require.

Fortunately, though, this slightly more complicated situation follows from Proposition 2.1. To see this, we note that $V(t)$ solves the Cauchy problem

$$\begin{cases}\left(\tfrac{1}{i}\tfrac{\partial}{\partial t} + P\right)V(t) = [P,C]U(t)B \\ V(0) = CB.\end{cases}$$

Consequently, by DuHamel's formula

(18) $$CU(t)B = U(t)CB + i\int_0^t U(t)\bigl(U(-s)[P,C]U(s)B\bigr)ds.$$

If we change our notation a bit and let $\sigma_{prin}(CB)$ and $\sigma_{sub}(CB)$ denote the principal and subprincipal symbols of $CB$ then Proposition 2.1 tells us that for small $|t|$ we can write the restriction to the diagonal of the kernel of $U(t)CB$ as

$$\int \frac{\partial A_0}{\partial\lambda}(x,\lambda)e^{-i\lambda t}\,d\lambda,$$

where

$$A_0(x,\lambda) - (2\pi)^{-n}\int_{p<\lambda}\bigl(\sigma_{prin}(CB) + \sigma_{sub}(CB)\bigr)d\xi$$
$$+ (2\pi)^{-n}\frac{\partial}{\partial\lambda}\int_{p<\lambda}\tfrac{i}{2}\{\sigma_{prin}(CB),p\}d\xi \in S^{n-2}.$$

To perform the same calculation for the last term in (18), we note that by Egorov's theorem
$$Q(t, x, D_x) = i \int_0^t U(-s)[P, C]U(s)B ds$$
is as in Proposition 2.1 with $Q(0, x, D_x) = 0$ and $\partial_t Q(0, x, D_x) = i[P, C]B$. Thus, for small $|t|$ the restriction to the diagonal of the kernel of the last term in (18) can be written as
$$\int \frac{\partial A_1}{\partial \lambda}(x, \lambda)e^{-i\lambda t} d\lambda,$$
where
$$A_1(x, \lambda) - (2\pi)^{-n} \frac{\partial}{\partial \lambda} \int_{p<\lambda} \sigma_{prin}([P, C]B) d\xi \in S^{n-2}.$$

Let $b$ and $c$ denote the principal symbols of $B$ and $C$. Note then that
$$\tfrac{i}{2}\{\sigma_{prin}(CB), p\} - \sigma_{prin}([P, C]B) = \tfrac{i}{2}\bigl(c\{b, p\} - b\{c, p\}\bigr).$$

Thus, we can combine the main term for $A_1$ with the last term for $A_0$ and obtain the following result which will play a key role in the proof of (11).

**Proposition 2.2.** *The restriction $K$ of the kernel of $CU(t)B$ to $\mathbb{R} \times \{(x, x) : x \in M\}$ is conormal with respect to $\{0\} \times \{(x, x) : x \in M\}$ in a neighborhood of this submanifold. Moreover, there is a $\delta > 0$ so that when $|t| < \delta$*

(19) $$K(t, x) = \int_{-\infty}^{\infty} \frac{\partial A_{CB}(x, \lambda)}{\partial \lambda} e^{-i\lambda t} d\lambda,$$

*where $A_{CB} \in S^n$, $A_{CB}(0, \lambda) = 0$ and*

(20) $$A_{CB}(x, \lambda) - (2\pi)^{-n} \int_{p(x,\xi)<\lambda} (cb + \sigma_{sub}(CB)) d\xi$$
$$+ (2\pi)^{-n} \frac{\partial}{\partial \lambda} \int_{p(x,\xi)<\lambda} \tfrac{i}{2}\bigl(c\{b, p\} - b\{c, p\}\bigr) d\xi \in S^{n-2}.$$

2.3. **Tauberian Lemmas.** In this section we shall collect the Tauberian lemmas that we need.

The first one is a special case of Lemma 17.5.6 in [Ho III]. It is a slight variant of the one used by Ivriĭ[Iv1] in his proof of the Duistermaat-Guillemin theorem. It requires a monotonicity assumption that will only be fulfilled for the "diagonal" terms of our approximation to $E_\lambda(x, x)$.

**Lemma 2.3.** *Suppose that $\mu$ is a non-decreasing temperate function satisfying $\mu(0) = 0$ and that $\nu$ is a function of locally bounded variation such that $\nu(0) = 0$. Suppose also that $m \geq 1$ and that $\phi \in \mathcal{S}(\mathbb{R})$ is a fixed positive function satisfying $\int \phi(\lambda)d\lambda = 1$ and $\hat{\phi}(t) = 0$, $t \notin [-1, 1]$. If $\phi_\sigma(\lambda) = \sigma^{-1}\phi(\lambda/\sigma)$, $0 < \sigma \leq \sigma_0$, assume that for $\lambda \in \mathbb{R}$*

(21) $$|d\nu(\lambda)| \leq \bigl(A_0(1+|\lambda|)^m + A_1(1+|\lambda|)^{m-1}\bigr) d\lambda,$$

*and that*

(22) $$|((d\mu - d\nu) * \phi_\sigma)(\lambda)| \leq B(1+|\lambda|)^{-2}.$$

*Then*

(23) $$|\mu(\lambda) - \nu(\lambda)| \leq C_m \big( A_0 \sigma(1+|\lambda|)^m + A_1 \sigma(1+|\lambda|)^{m-1} + B \big),$$

*where $C_m$ is a uniform constant depending only on $\sigma_0$ and our $m \geq 1$.*

The other lemma that we require will allow us to handle the "off-diagonal" terms in our approximation to $E_\lambda(x,x)$ where the above monotonicity assumption will not be valid. This lemma is exactly the one employed by Hörmander [Ho 1] and Avakumovič [Av] and Levitan [Le] in their proofs of the $O(\lambda^{n-1})$ remainder estimates for the Weyl formula. A proof can be found in [Ho 1] (or [So3], p. 128).

**Lemma 2.4.** *Let $g(\lambda)$ be a piecewise continuous function of $\mathbb{R}$. Assume that for some $m \geq 1$ there is a constants $A_0$ and $A_1$ so that*

(24) $$|g(\lambda + s) - g(\lambda)| \leq A_0(1+|\lambda|)^m + A_1(1+|\lambda|)^{m-1}, \quad 0 < s < 1.$$

*Suppose further that for some fixed $\delta > 0$*

$$\hat{g}(t) = 0, \quad |t| < \delta.$$

*Then there is a constant $C_{m,\delta}$, depending only on $m$ and $\delta$, so that*

(25) $$|g(\lambda)| \leq C_{m,\delta}\big(A_0(1+|\lambda|)^m + A_1(1+|\lambda|)^{m-1}\big).$$

We also quote a simple asymptotic result of [Z] (Lemma 3.6), which uses this Tauberian Lemma and which will be used in §7. In the following, $SV$ for any homogeneous $V \subset T_x^*M$ denotes the set of unit vectors in $V$. For background on Maslov indices, excesses etc. we refer to [DG] [Ho IV] [Z].

**Lemma 2.5.** *Let $A \in I^0(M \times M, \Lambda)$ denote a 0th order Fourier integral operator associated to a homogeneous canonical relation $\Lambda \subset T^*M \times T^*M$, and assume that $\Lambda \cap \Delta_{T^*M}$ is a clean intersection. Then*

$$\sum_{\lambda_\nu \leq \lambda} \langle A\phi_\nu, \phi_\nu \rangle = i^m \lambda^{\frac{e-1}{2}} \int_{S\Lambda_\Delta} \sigma_A d\mu_0 + O(\lambda^{\frac{e-1}{2}-1})$$

*where $d\mu_0$ is a canonical density on the intersection $\Lambda \cap \Delta_{T^*M}$, $\sigma_A$ is the symbol of $A$, $m$ is the Maslov index, and $e = \dim \Lambda \cap \Delta_{T^*M}$.*

## 3. Proofs of Theorems 1.1 and 1.2

We first concentrate on our Weyl formula remainder estimates in Theorem 1.2 and postpone proof of our eigenfunction bounds in Theorem 1.1 until the end of this section.

### 3.1. Local Weyl law remainder: Proof of Theorem 1.2. Following [Ho IV], we begin by introducing the loop-length function on $T^*M \backslash 0$ given by

(26) $$L^*(x,\xi) = \inf\{t > 0 : \exp_x t\xi = x\},$$

where $L^*$ is defined to be $+\infty$ if no such $t$ exists. It is homogeneous of degree zero, so it is natural to consider the restriction of $L^*$ to $S^*M = \{(x,\xi) : \sum g^{jk}(x)\xi_j\xi_k = 1\}$. A key fact for us is that $L^*$ is a lower semicontinuous function, or equivalently that the function $1/L^*(x,\xi)$, which is defined to be zero when $L^*(x,\xi) = +\infty$, is an upper semicontinuous function on $S^*M$. For fixed $x \in M$ we define the set of loop directions at $x$ by:

$$\mathcal{L}_x = \{\xi \in S^*_x M : 1/L^*(x,\xi) \neq 0\}. \tag{27}$$

The complimentary set

$$S^*M_x \backslash \mathcal{L}_x = \{\xi \in S^*_x M : 1/L^*(x,\xi) = 0\}$$

is the set of all unit vectors for which there is no geodesic loop with initial tangent vector $\xi$.

Let us first show that if $x_0 \in M$ and $|\mathcal{L}_{x_0}| = 0$, then given $\varepsilon > 0$, we can find a ball $B$ centered at $x_0$ and a $\Lambda < \infty$ so that

$$|R(\lambda,x)| \leq \varepsilon \lambda^{n-1}, \ x \in B, \ \lambda \geq \Lambda. \tag{28}$$

To do this, choose a coordinate patch $\Omega$, with coordinates $\kappa(x)$ around $x_0$ which we identify with an open subset of $\mathbb{R}^n$. Let then $f(x,\xi)$ denote the image of $1/L^*(x,\xi)$ in the induced coordinates for $\{(x,\xi) \in S^*M : x \in \Omega\}$. We then need the following

**Lemma 3.1.** *Let $f$ be a nonnegative upper semicontinuous function on $\mathcal{O} \times S^{n-1}$, where $\mathcal{O} \subset \mathbb{R}^n$ is open. Fix $x_0 \in \mathcal{O}$ and suppose that $\{\xi \in S^{n-1} : f(x_0,\xi) \neq 0\}$ has measure zero. Let $\varepsilon > 0$ and $\delta > 0$ be given. Then there is a neighborhood $\mathcal{N}$ of $x_0$ and an open set $\Omega_b$ satisfying*

$$|f(x,\xi)| < \delta, \quad (x,\xi) \in \mathcal{N} \times S^{n-1} \backslash \Omega_b$$
$$|\Omega_b| \leq \varepsilon.$$

*Furthermore, there is a $b(\xi) \in C^\infty$ supported in $\Omega_b$, satisfying $0 \leq b \leq 1$, and having the property that if $B(\xi) = 1 - b(\xi)$ then $|f(x,\xi)| < \delta$ on $\mathcal{N} \times \mathrm{supp}\ B$.*

**Proof.** By assumption,

$$E_\delta = \{\xi \in S^{n-1} : f(x_0,\xi) \geq \delta\}$$

satisfies $|E_\delta| = 0$. Let

$$E_\delta(j) = \{\xi \in S^{n-1} : f(x,\xi) \geq \delta \text{ some } x \in \overline{B}(x_0, 1/j)\},$$

if $\overline{B}(x_0, 1/j)\}$ is the closed ball of radius $1/j$ about $x_0$. Then clearly, $E_\delta(j+1) \subset E_\delta(j)$. Also, if $\xi \in \cap_{j \geq 1} E_\delta(j)$, then for every $j$ one can find $x_j \in B(x, 1/j)$ so that $f(x_j,\xi) \geq \delta$, which by the upper semicontinuity of $f$ in turn implies

$$\delta \leq \limsup_{j \to \infty} f(x_j,\xi) \leq f(x_0,\xi).$$

Thus,

$$\cap_{j \geq 1} E_\delta(j) \subset E_\delta.$$

Consequently, if $j$ is large, $|E_\delta(j)| < \varepsilon/2$. Fix such a $j = j_0$ and choose an open set $\Omega_b$ satisfying $E_\delta(j_0) \subset \Omega_b$ and $|\Omega_b| < \varepsilon$. Then clearly $|f(x,\xi)| < \delta$ if $(x,\xi) \in B(x_0, 1/j_0) \times S^{n-1} \backslash \Omega_b$.

For the last part, note that the argument we have just given will show that the sets $E_\delta(j)$ are closed because of the upper semicontinuity of $f$. Thus, if $E_\delta(j_0)$ and $\Omega_b$ are as above, we need only apply the $C^\infty$ Urysohn lemma to find a smooth function $b(\xi)$ supported in $\Omega_b$ with range $[0,1]$ and satisfying $b(\xi) = 1$, $\xi \in \Omega_b$, which then will clearly have the required properties. □

It follows from this Lemma that if $|\mathcal{L}_{x_0}| = 0$, then given $0 < \ell < \infty$ there cannot be too many geodesic loops of length $\leq \ell$ through points $x$ near $x_0$.

As in Ivriĭ's proof [Iv2] (see also [Ho IV]) of the Duistermaat-Guillemin theorem, we now construct a microlocal partition of unity which is adapted to the loop geometry of $(M, g)$. It will suffice to construct two pseudodifferential operators $b(x, D), B(x, D) = I - b(x, D)$ with the property that $B$ is microsupported in the set where $L^*(x, \xi) >> T$. We then study $E_\lambda(x, x)$ by writing it in the form:

$$(29) \qquad E_\lambda(x, x) = [(B + b)E_\lambda(B + b)](x, x).$$

The expert reader will observe that this expression is more complicated than the one

$$E_\lambda = B^* B E_\lambda + b^* b E_\lambda$$

used in [Iv1] [Ho IV]. However, the simpler form is not useful in our problem because we are not taking the trace of $E_\lambda$ and cannot move the adjoints in front of $E_\lambda$ to make the individual terms into positive operators. Rather, we shall have to expand (29) and study each time of term separately in Lemmas 3.2- 3.3. Moreover we shall have to define $b$ carefully to obtain the continuity properties of $R(\lambda, x)$. We now turn to the details of the construction.

If $T > 1$ is a large parameter, then by Lemma 3.1 we can find a function $b \in C^\infty(S^{n-1})$ with range $[0, 1]$ so that

$$\int_{S^{n-1}} b(\xi) d\sigma(\xi) \leq 1/T^2,$$

and

$$|f(x, \xi)| \leq 1/T, \quad \text{on } \mathcal{N} \times \text{supp } B,$$

where $\mathcal{N} \subset \kappa(\Omega)$ is a neighborhood of $\kappa(x_0)$ and

$$B(\xi) = 1 - b(\xi).$$

Choose a real-valued function $\psi \in C^\infty$ which vanishes outside of $\mathcal{N}$ and equals one in a small ball centered at $\kappa(x_0)$. Using these functions we get zero-order pseudo-differential operators on $\mathbb{R}^n$ by setting

$$\tilde{b}(x, D)f(x) = \psi(x)(2\pi)^{-n} \iint e^{i(x-y) \cdot \xi} b(\xi/|\xi|) \psi(y) f(y) \, d\xi dy,$$

and

$$\tilde{B}(x, D)f(x) = \psi(x)(2\pi)^{-n} \iint e^{i(x-y) \cdot \xi} B(\xi/|\xi|) \psi(y) f(y) \, d\xi dy.$$

Note that both variables of the kernels of these operators have compact support in $\mathcal{N}$. If we let $b(x, D)$ and $B(x, D)$ in $\Psi^0(M)$ be the pullbacks of $\tilde{b}(x, D)$ and $\tilde{B}(x, D)$, respectively, then

$$(30) \qquad b(x, D) + B(x, D) = \psi^2(x),$$

where we are abusing the notation a bit by letting $\psi$ denote multiplication by the pullback of our $\mathbb{R}^n$ cutoff function.

Note that by construction

$$(31) \qquad B(x, \xi) = 0 \text{ if } (x, x, t; \tau, \xi, \eta) \in \mathcal{C} \text{ for some } (t, \tau, \eta) \text{ with } 0 < |t| < T,$$

with $\mathcal{C}$ being the Lagrangean of $U(t,x,y)$. Consequently the kernels of $U(t)B^*$ and $BU(t)$ are smooth at the diagonal when $0 < |t| < T$.

Next, let $b_0$, $b^s$ and $B_0$, $B^s$ denote the principal and subprincipal symbols of $b(x,D)$ and $B(x,D)$, respectively. Then

$$(32) \qquad T^2 \int_{p(x,\xi)\leq 1} |b_0(x,\xi)|^2\,d\xi + \int_{p(x,\xi)\leq 1} |B_0(x,\xi)|^2\,d\xi \leq C,$$

where $C$ depends on our local coordinate transformation above but is independent of the large parameter $T$. Also, one has bounds

$$(33) \quad \int_{p(x,\xi)\leq 1} (|b_0| + |B_0|)(|\{b_0,p\}| + |\{B_0,p\}|)d\xi$$
$$+ \int_{p(x,\xi)\leq 1} \bigl(|\sigma_{sub}(BB^*)| + \sigma_{sub}(bb^*)| + \sigma_{sup}(Bb^*)|\bigr)d\xi\,d\xi \leq C_T,$$

where this constant does depend on $T$.

We shall use these facts and the results of the last two sections to prove that for $\lambda \geq 1$

$$(34) \qquad \bigl|\psi^4(x)\bigl(E_\lambda(x,x) - (2\pi)^{-n}\int_{p(x,\xi)\leq\lambda} d\xi\,\bigr)\bigr| \leq CT^{-1}\lambda^{n-1} + C_T\lambda^{n-2},$$

where $C_T$ depends on $T$ but $C$ does not. This of course yields (28).

In view of (30) this bound is a consequence of the following two lemmas.

**Lemma 3.2.** *There are constants $C$ and $C_T$ as in (34) so that if $B = B(x,D)$ and $B^*$ is its adjoint then for $\lambda \geq 1$*

$$(35) \quad |BE_\lambda B^*(x,x) - (2\pi)^{-n}\int_{p(x,\xi)\leq\lambda}(|B_0|^2 + \sigma_{sub}(BB^*))d\xi$$
$$+ (2\pi)^{-n}\partial_\lambda \int_{p(x,\xi)\leq\lambda} \tfrac{i}{2}(B_0\{B_0^*,p\} - B_0^*\{B_0,p\})d\xi|$$
$$\leq CT^{-1}\lambda^{n-1} + C_T\lambda^{n-2},$$

*and if $b = b(x,D)$ and $b^*$ is its adjoint then for $\lambda \geq 1$*

$$(36) \quad |bE_\lambda b^*(x,x) - (2\pi)^{-n}\int_{p(x,\xi)\leq\lambda}(|b_0|^2 + \sigma_{sub}(bb^*))d\xi$$
$$+ (2\pi)^{-n}\partial_\lambda \int_{p(x,\xi)\leq\lambda} \tfrac{i}{2}(b_0\{b_0^*,p\} - b_0^*\{b_0,p\})d\xi|$$
$$\leq CT^{-2}\lambda^{n-1} + C_T\lambda^{n-2}.$$

**Lemma 3.3.** *There are constants $C$ and $C_T$ as in (34) so that for $\lambda \geq 1$*

(37)
$$|BE_\lambda b^*(x,x) + bE_\lambda B^*(x,x) - (2\pi)^{-n} \int_{p(x,\xi)\leq\lambda} (B_0 b_0^* + b_0 B_0^* + \sigma_{sub}(Bb^*) + \sigma_{sub}(bB^*))d\xi$$
$$+ (2\pi)^{-n} \partial_\lambda \int_{p(x,\xi)\leq\lambda} \tfrac{i}{2}(B_0\{b_0^*, p\} - b_0^* B_0, p + b_0\{B_0^*, p\} - B_0^*\{b_0, p\})d\xi|$$
$$\leq CT^{-1}\lambda^{n-1} + C_T \lambda^{n-2}.$$

One obtains (34) from these three inequalities after noting that (30) implies that
$$(B_0 + b_0)(B_0^* + b_0^*) = \psi^4,$$
and
$$\sigma_{sub}(BB^*) + \sigma_{sub}(bb^*) + \sigma_{sub}(Bb^*) + \sigma_{sub}(b^*B) = \sigma_{sub}(\psi^4) = 0.$$

**Proof of Lemma 3.2** To prove (35) we note that by (31) $K(t,x) = BU(t)B^*(x,x)$ is smooth for $0 < |t| < T$. Therefore if we let $\nu(\lambda) = A_{BB^*}(x,\lambda)$ be as in (20), then since (19) holds for small values of $t$ if $\mu(\lambda) = BE_\lambda B^*(x,x)$ then we have
$$|((d\mu - d\nu) * \phi_\sigma)(\lambda)| \leq B_N(1+|\lambda|)^{-N}$$
for constants $B_N$ independent of $N$ if $\sigma \leq 1/(2T)$ is fixed. By (32) and (33), (21) must hold with $A_0 = C$ and $A_1 = C_T$ where $C$ and $C_T$ are as in (34). Since
$$\lambda \to \mu(\lambda) = \sum_{\lambda_\nu \leq \lambda} |B\phi_\nu|^2$$
is non-decreasing Lemma 2.3 yields (35).

The proof of (36) follows from the same lines except that here one takes $\sigma$ to be twice the injectivity radius of $M$ and uses the fact that $\int_{p<1} |b_0|^2 d\xi \leq C/T^2$ to get the $O(T^{-2}\lambda^{n-1}) + O(\lambda^{n-2})$ bounds. □

**Proof of Lemma 3.3** In this case since
$$\lambda \to BE_\lambda b^*(x,x) + bE_\lambda B^*(x,x)$$
is no longer non-decreasing the first Tauberian lemma, Lemma 2.3, does not apply. Instead we must use the second one, Lemma 2.4.

To this end let $g(\lambda) = g_x(\lambda)$ be given by
$$g(\lambda) = \big(BE_\lambda b^*(x,x) + bE_\lambda B^*(x,x)\big) - \big(A_{Bb}^*(x,\lambda) + A_{bB^*}(x,\lambda)\big),$$
where $A_{Bb^*}$ and $A_{bB^*}$ are given by (20). By Proposition 2.2, the Fourier transform of $dg$ vanishes for $|t| < \delta$, and so supp $t\hat{g}(t) \subset \{0\} \cup \{t : |t| > \delta\}$. But then if $\{0\} \in$ supp $\hat{g}$, $g$ would include a sum of derivatives of delta functions, which is impossible since $g$ is clearly a piecewise continuous function, and so $\hat{g}(t) = 0$ for $|t| < \delta$. Therefore, by Lemma 2.4 we would have (37) if we could show that there are constants $C$ and $C_T$ as in (34) so that

(38) $\quad |g(\lambda+s) - g(\lambda)| \leq CT^{-1}(1+|\lambda|)^{n-1} + C_T(1+|\lambda|)^{n-2}, \ 0 < s \leq 1.$

Let us first see that this bound holds when $g$ is replaced by $A_{Bb^*}$ or $A_{bB^*}$. Here, we shall just use the fact that, modulo terms of order $n-2$, both $A_{Bb^*}$ and $A_{bB^*}$ are the sum of a function which is homogeneous of degree $n$ and one of degree $n-1$. Consequently, the bound will just follow from the size of the coefficients of these two homogeneous terms. For $A_{Bb^*}$ we note that, by (20),

$$A_{Bb^*}(x,\lambda) = c(x)\lambda^n + d(x)\lambda^{n-1} + O(\lambda^{n-2}),$$

where

$$c(x) = (2\pi)^{-n}\int_{p(x,\xi)\leq 1} B_0 \overline{b_0}\, d\xi,$$

and

$$d(x) = (2\pi)^{-n}\int_{p(x,\xi)\leq 1} \sigma_{sub}(Bb^*)d\xi - n(2\pi)^{-n}\int_{p(x,\xi)\leq 1} \tfrac{i}{2}\big(B_0\{\overline{b}_0,p\} - \overline{b}_0\{B_0,p\}\big)d\xi.$$

By (32) and Schwarz's inequality, $|c(x)| \leq C/T$, while (33) yields $|d(x)| \leq C_T$ and so (38) must hold when $g$ is replaced by $A_{Bb^*}$.

Clearly the same argument applies to $A_{bB^*}$. Therefore, since

$$|BE_{\lambda+s}b^*(x,x) - BE_\lambda b^*(x,x)| = |bE_{\lambda+s}B^*(x,x) - bE_\lambda B^*(x,x)|$$
$$= |\sum_{\lambda<\lambda_\nu\leq\lambda+s} B\phi_\nu(x)\overline{b\phi_\nu(x)}|,$$

the proof would be complete if we could show that

$$|\sum_{\lambda<\lambda_\nu\leq\lambda+s} B\phi_\nu(x)\overline{b\phi_\nu(x)}| \leq CT^{-1}(1+|\lambda|)^{n-1} + C_T(1+|\lambda|)^{n-2},\ 0\leq s\leq 1.$$

To prove this we apply the Cauchy-Schwarz inequality to get

$$|\sum_{\lambda<\lambda_\nu\leq\lambda+s} B\phi_\nu(x)\overline{b\phi_\nu(x)}|$$
$$\leq (\sum_{\lambda<\lambda_\nu\leq\lambda+1}|B\phi_\nu|^2)^{1/2}(\sum_{\lambda<\lambda_\nu\leq\lambda+1}|b\phi_\nu(x)|^2)^{1/2}$$
$$= |BE_{\lambda+1}B^*(x,x) - BE_\lambda B^*(x,x)|^{1/2}\,|bE_{\lambda+1}b^*(x,x) - bE_\lambda b^*(x,x)|^{1/2}.$$

By (35) and the argument just given for $A_{Bb^*}$

$$|BE_{\lambda+1}B^*(x,x) - BE_\lambda B^*(x,x)| \leq C(1+|\lambda|)^{n-1} + C_T(1+|\lambda|)^{n-2},$$

since the coefficient of the term in the left side of (35) which is homogeneous of degree $n$ is bounded by a uniform constant $C$, while that of the $(n-1)$ degree term is bounded by a constant $C_T$.

To finish we need to see therefore that

$$|bE_{\lambda+1}b^*(x,x) - bE_\lambda b^*(x,x)| \leq CT^{-2}(1+|\lambda|)^{n-1} + C_T(1+|\lambda|)^{n-2}.$$

but this follows from (36) since the coefficient of the term which is homogeneous of degree $n$ is bounded by $CT^{-2}$.

This completes the proof of Theorem 1.2.

3.2. **Applications to eigenfunctions: Proof of Theorem 1.1.** It is now a simple matter to complete the proof of Theorem 1.1. We begin with the following simple fact.

**Lemma 3.4.** *Fix $x \in M$. Then if $\lambda \in \text{spec } \sqrt{-\Delta}$*

$$\sup_{\phi \in V_\lambda} \frac{|\phi(x)|}{\|\phi\|_2} = \sqrt{R(\lambda, x) - R(\lambda - 0, x)}, \qquad (39)$$

*if $R(\lambda - 0, x)$ denotes the left-hand limit of $R$ at $\lambda$.*

**Proof.** Let
$$K_\lambda(x, y) = \sum_{\lambda_\nu = \lambda} \phi_\nu(x) \overline{\phi_\nu(y)}$$
denote the kernel of the projection onto $V_\lambda$. Then if $x \in M$ is fixed, $y \to K_\lambda(x, y)$ is in $V_\lambda$. Consequently, by orthogonality,

$$\Big( \sum_{\lambda_\nu = \lambda} |\phi_\nu(x)|^2 \Big)^{1/2} = \Big( \int |K_\lambda(x, y)|^2 \, dy \Big)^{1/2} = \sup_{\phi \in V_\lambda, \|\phi\|_2 = 1} \Big| \int K_\lambda(x, y) \phi(y) \, dy \Big|$$
$$= \sup_{\phi \in V_\lambda, \|\phi\|_2 = 1} |\phi(x)|.$$

Note that $\lambda \to E_\lambda(x, x)$ is an increasing right continuous function, and so by (4)
$$\sum_{\lambda_\mu = \lambda} |\phi_\nu(x)|^2 = R(\lambda, x) - R(\lambda - 0, x),$$
which along with the last identity gives (39). $\square$

The first bound in the statement of Theorem 1.1 is then an immediate consequence of (12) and (39). To prove the remaining estimate, (8), let us note that it is equivalent to showing that if $|\mathcal{L}_x| = 0$ for every $x \in M$, then if $\varepsilon > 0$ there is a $\Lambda(\varepsilon, p)$ so that

$$\|e_\lambda(f)\|_p \le \varepsilon \lambda^{\delta(p)} \|f\|_2, \quad \lambda \ge \Lambda(\varepsilon, p), \quad \tfrac{2(n+1)}{n-1} < p \le \infty, \qquad (40)$$

if $e_\lambda$ is the projection onto $V_\lambda$. The endpoint $p = \infty$ follows from (12) and Lemma 3.1. Since it was shown in [So1] that one has the uniform $O$-bounds

$$\|e_\lambda(f)\|_p \le C \lambda^{\delta(p)} \|f\|_p, \quad p = \tfrac{2(n+1)}{n-1}, \qquad (41)$$

if we apply the M. Riesz interpolation theorem (see, e.g. [So3] p. 6) to the $p = \infty$ special case of (40) and (41) we get

$$\|e_\lambda(f)\|_p \le C^{\theta_p} \varepsilon^{1-\theta_p} \lambda^{\delta(p)} \|f\|_p, \quad \lambda \ge \Lambda(\varepsilon, \infty), \ p > \tfrac{2(n+1)}{n-1}, \ \theta_p = \tfrac{2(n+1)}{(n-1)p},$$

which of course yields the remaining cases of (41). $\square$

3.3. **Duistermaat-Guillemin.** Since it is quite simple, for the sake of completeness, let us point out how our pointwise bounds for the Weyl remainder term imply integrated ones. Specifically, if $N(\lambda)$ is the number of eigenvalues of $\sqrt{-\Delta}$ that are $\le \lambda$, then our results imply that

$$N(\lambda) = (2\pi)^{-n} \text{Vol } M \lambda^n + R(\lambda), \quad R(\lambda) = o(\lambda^{n-1}),$$

under the assumption that $|\mathcal{L}| = 0$, if $\mathcal{L} \subset S^*M$ is the set of initial directions $(x, \xi)$ of closed loops.

To see this we note that our assumption $|\mathcal{L}| = 0$ and the Tonelli-Fubini theorem imply that one must have $|\mathcal{L}_x| = 0$ when $x \in M \backslash E$ with $|E| = 0$. To use this, if $\varepsilon > 0$ is fixed, choose an open set $\Omega \subset M$ containing $E$ with $|\Omega| < \varepsilon$. Since $K = M \backslash \Omega$ is compact, we conclude from (12) that for this $\varepsilon$ we can find a $\Lambda \geq 1$ so that $|R(\lambda, x)| \leq \varepsilon \lambda^{n-1}$ when $x \in K$ and $\lambda \geq \Lambda$. On other hand, we know from [Ho 1] that there is a uniform constant $C$ so that $|R(\lambda, x)| \leq C\lambda^{n-1}$, $\lambda \geq 1, x \in M$. Thus, if $\lambda \geq \Lambda$,

$$|R(\lambda)| = |\int_M R(\lambda, x)dx| \leq \int_\Omega |R(\lambda, x)|dx + \int_K |R(\lambda, x)|dx$$
$$\leq C\varepsilon \text{Vol } M\lambda^{n-1} + \varepsilon \text{Vol } M\lambda^{n-1},$$

which of course means that $R(\lambda) = o(\lambda^{n-1})$, as desired.

We should point out that the assumption $|\mathcal{L}| = 0$ seems to be stronger than the assumption of Duiestermaat-Guillemin [DG] that the set of closed geodesics has measure zero. The two assumptions might be the same, but we do not know of a proof.

## 4. Geometry of loops

In the remainder of the paper, we apply Theorem 1.1 to a number of concrete classes of compact Riemannian manifolds and investigate the extent to which it is sharp. In so doing, we shall need to analyze the condition $|\mathcal{L}_x| = 0$ ($\forall x \in M$) and related properties of loops. We therefore pause in our analysis of eigenfunctions to explore the geometry of loops. The issues involved are motivated by analytic problems and lie someone outside of the mainstream of the geometry of geodesics.

To clarify the geometry of loops, we introduce the following notions:

**Definition 4.1.** *We put:*

- $\Gamma_x^T := \exp T H_p T_x^* M \cap T_x^* M$.
- $S\Gamma_x^T = \Gamma_x^T \cap S_x^* M$.
- $\Gamma_x = \exp_x^{-1}(x) = \{T\xi : \xi \in S\Gamma_x^T\} \subset T_x^* M$;
- $Lsp_x = \{|T| : \Gamma_x^T \neq \emptyset\} \subset \mathbb{R}^+$
- $\mathcal{L}_x := \{\xi \in S_x^* M : \exists t \in Lsp_x, \exp_x t\xi = x\} \subset S_x^* M$

In this list, we have mixed together notions involving the geometer's geodesic flow (the Hamiltonian flow of $p^2$) and the microlocal analyst's geodesic flow (that of $p$). Thus, $\Gamma_x^T$ is a homogeneous closed subset of $T_x^* M$ (invariant under the $\mathbb{R}^+$-action), and $S\Gamma_x^T$ is the 'base' of the cone. On the other hand, we define $\exp_x$ in the usual geometer's, rather than as $\pi_x \circ \exp tH_p|_{T_x^* M}$. For instance, we have: $TS\Gamma_x^T = \{T\xi : |\xi| = 1, \exp xT\xi = x\}$.

It is natural to say that $Lsp_x$ is the 'loop-length spectrum' at the point $x$, i.e. the set of lengths of geodesic loops. The set $\mathcal{L}_x$, already introduced, is the set of loop directions at $x$. The (non-homogeneous) set $\Gamma_x$ is the set of loop (co)-vectors. In general, the sets $\Gamma_x^t$, $\Gamma_x$ need not be smooth manifolds. Later in this section, we shall consider restrictions on $(M, g)$ which imply smoothness properties of loopsets, and study a variety of examples.

First, we present a standard result which plays a crucial role in several applications.

**Proposition 4.2.** $L^*$ *is constant on smooth submanifolds of* $\Gamma_x$.

*Proof.* It is sufficient to prove the result for smooth curves in $\Gamma_x$. This can be seen by using either the first variation formula for the lengths of a one-parameter family of geodesics, or by using a symplectic argument. We present both for the reader's convenience.

**First Variation proof**: Let $\alpha(s)$ denote a smooth curve in $\Gamma_x$, and let $\gamma_s(t) = \exp_x tL^*(x, \alpha(s)) \frac{\alpha(s)}{|\alpha(s)|}$. By definition, $\gamma_s(0) = x, \gamma(1) = x$ and $L^*(x, \alpha(s)) = \int_0^1 |\gamma'_s(t)| dt$. Since $\gamma_s$ is a geodesic for each $s$, we have

$$\frac{d}{ds} \int_0^1 |\gamma'_s(t)| dt = \frac{1}{L^*(x,\alpha(s))} \int_0^1 \langle \frac{D}{ds} \gamma'_s(t),\ \gamma'_s(t) \rangle dt$$

$$\frac{1}{L^*(x,\alpha(s))} \int_0^1 \langle \frac{D}{dt} Y_s(t),\ \gamma'_s(t) \rangle dt$$

$$= -\frac{1}{L^*(x,\alpha(s))} \int_0^1 \langle Y_s(t),\ \frac{D}{dt} \gamma'_s(t) \rangle dt = 0.$$

Here, $Y_s(t) = \frac{d}{ds} \gamma_s(t)$ is the Jacobi field associated to the variation, and $\frac{D}{Dt}$ is covariant differentiation along $\gamma_s$. In integrating $\frac{D}{dt}$ by parts, the boundary terms vanished because $Y_s(0) = Y_s(1) = 0$ as the endpoints of the variation were fixed. □

**Symplectic proof**: The curve $C$ parameterized by $\alpha_s$ is an isotropic submanifold of $T^*M$, since it lies in $T_m^*M$. Hence, the map

$$i : C \to T^*M, \quad i(t\xi) = \exp(tH_p(m, \xi))$$

is an isotropic embedding. Now let $\alpha = \xi \cdot dx$ denote the action form of $T^*M$. Then the length of the geodesic loop $\gamma_s(t)$ (as above) is given by

$$L^*(x, \alpha(s)) = \int_{\gamma_s} \xi \cdot dx.$$

Let $S$ denote the surface parameterized by $\alpha(s,t) = \gamma_s(t)$. Since all the geodesic rays are loops, the boundary of $S$ equals $\gamma_{s_1} - \gamma_{s_0}$. Since $S$ is isotropic we have, by Stokes' theorem,

(42) $$L^*(x, \alpha(s_1)) - L^*(x, \alpha(s_0)) = \int_{\gamma_{\xi s_1}} \alpha - \int_{\gamma_{s_0}} \xi \cdot dx = \int_S \omega = 0.$$

Here, $\omega$ is the canonical symplectic form $dx \wedge d\xi$. □

**Corollary 4.3.** *Suppose that $\Gamma_x$ is a manifold. If $|\mathcal{L}_x| > 0$, then there exists $T > 0$ such that $TS_x^*M \subset \Gamma_x^T$.*

*Proof.* We have just proved that the connected components of $\Gamma_x$ must lie in the spheres $TS_x^*M$ where $\ell \in Lsp_x$. It also follows that the loop-length spectrum $Lsp_x$ is discrete, hence countable. Now, $\mathcal{L}_x$ is just the radial projection of $\Gamma_x$. If it has positive measure, then one of the components must project to a set of positive measure, hence $\Gamma_x^T$ must have positive measure in $TS_x^*M$. Since it is also a submanifold of $TS_x^*M$, it must equal it. □

We now introduce hypotheses on $(M, g)$ under which loopsets have good smoothness properties. We recall that the intersection $X \cap Y$ of two manifolds is *clean* if $X \cap Y$ is a manifold and if $T_x(X \cap Y) = T_x X \cap T_x Y$ for all $x \in X \cap Y$.

**Definition 4.4.** *We shall say that a compact Riemannian manifold $(M, g)$ has:*

(i) *Clean loops of length $T \neq 0$ if $\Gamma_x^T = \exp T H_p T_x^* M \cap T_x^* M$ is a clean intersection.*

(ii) *Clean loopsets at $x$ if $\Gamma_x^T$ is clean for all $T$.*

Cleanliness is related to the question of whether a one parameter family of geodesics at a self-conjugate point $x$ is a family of loops at $x$ or not. Let us note that a curve in $\Gamma_x^T$ is the set of initial vectors to a one-parameter family of geodesic loops $\gamma_u(t)$ ($t \in [0, T]$). The variation vector field $Y(t) = \frac{d}{du}\big|_{u=0} \gamma_u'(t)$ of such a family is then a Jacobi field satisfying $Y(0) = Y(T) = 0$. Thus, each tangent vector $V = \frac{d}{du}\gamma_u'(0)$ to $\Gamma_x^T$ has the form $V = Y'(0)$ where $Y$ is a Jacobi field satisfying $Y(0) = Y(T) = 0$. On the other hand, any Jacobi field $Y$ along a geodesic arc $\exp_x t\xi, t \leq T$ satisfying $Y(0) = Y(T) = 0$ gives rise to a one-parameter family of geodesics $\gamma_u(t), t \in [0, T]$ along which $x$ is self-conjugate. However, in general $\gamma_u(t), t \in [0, T]$ fails to be a loop at $x$ and is just a loop to 'first order'. Under the cleanliness condition, this cannot happen.

**Proposition 4.5.** $\exp T H_p T_x^* M \cap T_x^* M$ *is a clean intersection if and only if $\Gamma_x^T$ is a manifold and $T_{(x,\xi)}\Gamma_x^T$ is isomorphic to the space $\mathcal{J}_T$ of Jacobi fields $Y(t)$ along the geodesic arc $\exp_x t\xi, t \leq T$ satisfying $Y(0) = Y(T) = 0$.*

*Proof.* We obviously have
$$T_{(x,\xi)} \exp T H_p T_x^* M = (d_\eta \exp T H_p) T_{(x,\eta)} T_x^* M, \quad (\exp T H_p(x,\eta) = (x,\xi).$$
As is well known, $(d_\eta \exp T H_p)$ is defined by
$$(d_\eta \exp T H_p)(J(0), J'(0)) = (J(T), J'(T))$$
where $J(t)$ is a Jacobi field with initial conditions $(J(0), J'(0))$ along $\exp t H_p(x, \xi), t \in [0, T]$. Clearly, $\exp T H_p T_x^* M \cap T_{(x,\xi)} T_x^* M$ is the space of value $(J(T), J'(T))$ of endpoint values of Jacobi fields satisfying $J(0) = 0$. It follows that $T_{(x,\xi)} \exp T H_p T_x^* M \cap T_{(x,\xi)} T_x^* M$ consists of the Jacobi fields along the loop $\exp t H_p(x, \xi), t \in [0, T]$ with $J(0) = 0 = J(T)$. The cleanliness assumption is therefore $T_{(x,\xi)}\Gamma_x^T$ consists of all such Jacobi fields. □

We may reformulate the condition in terms of the length functional
$$L_x : \Omega_x \to \mathbb{R}, \quad \Omega = \{\alpha : [0,1] \to M, \alpha(0) = \alpha(1) = x, \alpha \in H^1[0,T]\}$$
on the $H^1$-loop space at $x$. As is well-known, the critical points of $L_x$ are just the loops at $x$, and the set of $T$ for which $\Gamma_x^T \neq \emptyset$ are the critical values of $L_x$.

**Proposition 4.6.** $(M, g)$ *has clean loopsets at $x$ if and only if $L_x$ is a Bott-Morse function on $\Omega_x$ i.e. a function whose critical point sets are non-degenerate critical manifolds.*

*Proof.* As noted in Proposition 4.2, $L_x$ is constant on each component of its critical point set, so the critical point set of critical value $T$ is precisely $\Gamma_x^T$. To prove non-degeneracy, we must show that the Hessian of $L_x$ is non-degenerate in the normal bundle of $\Gamma_x^T$ in $T_x^* M$. The Hessian at a critical point $xi$ is of course the Jacobi operator $D^2 + R(T, \cdot)T$ along the geodesic through $(x, \xi)$, restricted to tangent vectors to $\Omega_x$, i.e. to vector fields on the geodesic which vanish at both endpoints. Non-degeneracy of the Hessian is equivalent to the statement that its kernel lie in the tangent space to $\Gamma_x^T$, which we have just verified. □

In the general case (without cleanliness assumptions), it follows from Sard's theorem (and a standard finite dimensional approximation to $\Omega_x$ [M]) that $Lsp_x$ is always a set of measure zero in $\mathbb{R}$. Putting together the previous propositions, we obtain a stronger statement in the clean case.

**Corollary 4.7.** *Suppose that $(M,g)$ has clean loopsets. Then $Lsp_x$ is a discrete subset of $\mathbb{R}$ and $\Gamma_x$ is a submanifold of $T_x^*M$.*

*Proof.* Since $\Gamma_x = \cup_{T>0} \Gamma_x^T$ the second statement follows from the first and the definition of cleanliness. The first statement follows from the fact that Bott-Morse functions have isolated critical values. □

The following Proposition will play a role in the proof of Theorem 1.5:

**Proposition 4.8.** *Let $(M,g)$ be a compact Riemannian manifold. Then:*

*(i) Suppose that for some $T > 0$, the loopset $S\Gamma_x^T$ is clean and $|S\Gamma_x^T| > 0$. Then $S_x^*M = S\Gamma_x^T$.*

*(ii) If $(M,g)$ has clean loopsets at $x$, and if $|\mathcal{L}_x| > 0$, then there exists $T > 0$ such that $S_x^*M = S\Gamma_x^T$.*

*Proof.* (i) The hypothesis implies that $S\Gamma_x^T$ is a submanifold of $S_x^*M$ of full dimension.

(ii) Since $Lsp_x$ is discrete, $\Gamma_x = \cup_{T>0} T\Gamma_x^T$. Since $\mathcal{L}_x$ is the countable union of their radial projections to $S_x^*M$, and has positive measure, at least one component $T\Gamma_x^T$ must be of dimension $n-1$. Since $T\Gamma_x^T \subset TS_x^*M$ by Proposition 4.2, again it must equal it. □

**4.1. Examples.** We now present a few simple examples to illustrate the loopset notions, mainly with a view of clarifying the notation and terminology. It would take us too far afield to provide even a sketchy discussion of the wide variety of possible phenomena that occur; we limit ourselves to some rather standard cases.

4.1.1. *Flat tori.* Let $\mathbb{R}^n/\Gamma$ be a flat torus. A geodesic loop at $x$ is a helix which returns to $x$. Loops on flat tori are always closed geodesics and they correspond to lattice points in $\Gamma$. Thus,
$$\Gamma_x = \Gamma, \quad \mathcal{L}_x = \{\frac{\gamma}{|\gamma|} : \gamma \in \Gamma\}, \quad Lsp_x = \{|\gamma|, \ \gamma \in \Gamma\}, \quad \}.$$
We observe that $\Gamma_x$ is a countable discrete subset of $T_x^*\mathbb{R}^n/\Gamma$, that $\mathcal{L}_x$ is a countable dense subset of $S_x^*M$, and that $\Gamma_g$ is a countable set of embedded Lagrangean components each diffeomorphic to $\mathbb{R}^n/\Gamma$. Also, $\mathcal{C}_\Delta = \cup_{\gamma \in \Gamma}\{(|\gamma|, \tau, x, 0) : x \in M\}$. This example is non-generic because (among other things) every loop is a closed geodesic. Flat tori are examples of metrics without conjugate points. We now generalize the discussion.

4.1.2. *Manifolds without conjugate points.* Let $(M,g)$ denote a manifold without conjugate points, i.e. a Riemannian manifold such that each exponential map $\exp_x$ is is a covering map. Manifolds of non-positive curvature are examples, so this class of metrics is open on any $M$ (thought it may be empty for some $M$). By definition, there are no Jacobi fields along any loop satisfying $Y(0) = Y(1) = 0$, so the Jacobi operator is non-degenerate. Equivalently, each loop functional $L_x$ is a Morse function on $\Omega_x$. Thus, for all $x$, $\Gamma_x$ is a countable discrete set of points, $\mathcal{L}_x$ is countable and $\Gamma_x^T$ is a finite set. Unlike the case of flat tori, loops are not generally closed geodesics.

4.1.3. *Surfaces of revolution.* Surfaces of revolution provide a simple (but non-trivial) class of examples for which the loopsets may be explictly determined. Topologically, the surfaces must be either spheres or tori. It would take us too far afield to discuss the geometry of loops and eigenfunctions on general surfaces of revolution, so we only give a few indications of how to determine the loops explicitly and refer the reader to the articles [CV] [Bl][KMS][TZ] [Z2] for further background.

A 2-sphere of revolution is a Riemannian sphere $(S^2, g)$ with an action of $S^1$ by isometries. We may write $g = dr^2 + a(r)^2 d\theta^2$ in geodesic polar coordinates at one of the poles (fixed points) of the $S^1$ action. The poles are always self-conjugate points, and are points at which eigenfunctions attain the maximal bounds.

Tori of revolution generalize flat tori but non-flat cases must have conjugate points (Hopf's theorem). The metric on a torus of revolution may be written in standard angle coordiantes as $g = dx^2 + a(x)^2 d\theta^2$. Tori of revolution have no poles, and as we will see in Theorem 1.6 eigenfunctions need not attain the maximal bounds. In fact, although we do not prove it here, tori of revolution *never* have maximal eigenfunction growth.

The key feature of surfaces of revolution is the integrability of the geodesic flow, which follows from the existence of the Clairaut integral $I(x,\xi) = \xi(\frac{\partial}{\partial \theta})$. Here, $\frac{\partial}{\partial \theta}$ is the vector field generating the rotations. The Clairaut integral is constant along geodesics, or otherwise put, unit speed geodesics must lie on level sets of $I$ within $S_g^* M$. Compact regular level sets are necessarily tori and are known as the 'invariant tori.' We denote by $T_{I_0} \subset S_g^* M$ a torus with Clairaut integral $I = I_0$. For notational simplicity, we ignore the fact that the level sets of $I$ may have several components. Near any regular torus $T_I$, one can always find local 'action-angle' such that the geodesics may be expressed as winding lines on the invariant tori $(\phi_1 + t\omega_1, \phi_2 + t\omega_2)$. Here, $\omega_j$ are the frequencies of motion (which generally vary with $I$).

Now consider the geodesic loops at a point $x \in S^2$. If $x$ is a pole (a fixed point of the rotation), then all the geodesics passing through $x$ are meridians, and the projection from the torus of meridians $\{I = 0\}$ to $S^2$ under the standard projection $\pi : T^*S^2 \to S^2$ collapses the initial tangent vectors at $x$ to the meridians to the point $x$. Thus, $\Gamma_x^{2\pi} = T_x^* M$, $\mathcal{L}_x = S_x^* S^2$, and $Lsp_x = \{nL, n \in \mathbb{Z}\}$ where $L$ is the length of a closed meridian geodesic. All of these loops are closed geodesics.

Now let $M = S^2$ or $M = T^2$ and assume $x$ is not a pole, as must occur in the case of tori of revolution. Under $\pi : T^*M \to M$, the tori $T_I$ project to annuli on $M$, and with a finite number of exceptions, the projection map $\pi_I : T_I \to M$, restricted to a torus, is a 2-1 cover with fold singularities (the caustic sets). It follows there are either zero, one or two geodesics of $T_I$ passing through $x$, accordinly as $x$ does not lie in $\pi_I(T_I)$, lies on its boundary, or lies in its interior.

Using action-angle variables, we can determine the tori $T_I$ which contain a geodesic loop $x$. Such a loop is a winding line on $T_I$ which passes through the set $\pi_I^{-1}(x)$ at least twice. If it passes through the same point of $\pi_I^{-1}(x)$ twice, then it is a closed geodesic; otherwise it is a non-smooth loop. The set of tori which contain periodic geodesics are known as periodic tori, and are independent of the point $x$. Under a generic twist condition (see [Bl] [Z2]), the set of periodic tori in $S_g^*M$ is countable. It follows that for generic surfaces of revolution, the set of closed geodesics through any point is countable. More complicated are loops which are not closed geodesics, i.e. winding lines which pass through distinct points of $\pi_I^{-1}(x)$. These do depend on $x$. One can write down explicit equations for the values of $I$ such that $T_I$ contains a winding line which loops at a given point $x$. It is clear from the equations that, for fixed non-polar $x$, a generic surface of revolution has only countably many loops at $x$. However, the question whether $\Gamma_x$ is countable at all non-polar $x$ for generic tori of revolution appears to be open. For further information on geodesics on tori of revolution we refer to [KMS]. We plan to discuss this example in more detail in the future (both from the geometric and analytic point of view).

4.1.4. *Tri-axial ellipsoids.* We recall that $E_{a_1,a_2,a_3} = \{(x_1, x_2, x_3) \in \mathbb{R}^3 : \frac{x_1^2}{a_1^2} + \frac{x_2^2}{a_2^2} + \frac{x_3^2}{a_3^2} = 1\}$ with $0 < a_1 < a_2 < a_3$. Jacobi proved that the geodesic flow of $E_{a_1,a_2,a_3}$ (for any $(a_1, a_2, a_3)$) is completely integrable in 1838. The two integrals of the motion are the length $H(x, \xi) = |\xi|_x$ and the so-called Joachimsthal integral $J$. More recent discussions of the geodesics of the ellipsoid can be found in [A], [K], and in [CVV] (§3).

There are four distinguished umbilic points $\pm P, \pm Q$ which occur on the middle closed geodesic $\{x_2 = 0\}$. All geodesics leaving $P$ arrive at $-P$ at the same time, then leave $-P$ and return to $P$ at the same time (see [K], Theorem 3.5.16 or [A]). Thus, the tri-axial ellipsoid is an example of a $Y_\ell^m$-metric which is not a Zoll metric. At all other points $x \in E_{a_1,a_2,a_3}$, the set of initial directions of geodesics which return to $x$ is countable [K].

The moment map $(H, J) : T^*E_{a_1,a_2,a_3} \to \mathbb{R}^2$ is regular away the middle closed geodesic, and its only singular level in $S^*E_{a_1,a_2,a_3}$ is that $(H = 1, J = 0)$ of $\{x_2 = 0\}$. For a description of the singular level we refer to [CVV]. The middle geodesic is homoclinic, i.e. the trajectories on the level are forward and backward asymptotic to it. Other level sets of the moment map are regular Lagrangean tori, and the discussion of periodic orbits and loops is similar to the case of surfaces of revolution.

4.1.5. *Zoll surfaces.* Now let us consider the extreme case of Zoll metrics on $S^2$, i.e metrics all of whose geodesic are closed [Besse]. Among such metrics, there is an infinite dimensional family of surfaces of revolution. There is an even larger class with no isometries.

We may suppose with no loss of generality that the least common period of the geodesics equals $2\pi$. Then $\Gamma_x^{2\pi} = S_x^*S^2$ for every $x$, and $\Gamma_x = \cup_{n=1}^\infty 2\pi n S_x^*S^2$. As $x$ varies we obviously get $\Gamma = \cup_{n=1}^\infty 2\pi n S^*S^2$. We also note that $\mathcal{C}_\Delta$ is parameterized by

$$\mathbb{R} \times T^*M \to T^*(\mathbb{R} \times M \times M), \quad (t, x, \xi) \to (t, -|\xi|_x, x, \xi, \exp -tH_p(x, \xi)).$$

Since all geodesics are periodic, we obtain

$$\mathcal{C}_\Delta = \{(t, \tau, x, 0)\} \equiv \mathbb{Z} \times \mathbb{R} \times M,$$

clearly a homogeneous Lagrangean submanifold of $T^*(\mathbb{R} \times M)\backslash 0$.

As will be discussed in §8, although geodesics are recurrent at every $x$, we do not expect eigenfunction blow-up to occur everywhere, or even anywhere in general.

5. REAL ANALYTIC RIEMANNIAN MANIFOLDS

Our first application is to characterize Riemannian manifolds in the analytic setting with maximal eigenfunction growth.

THEOREM **5.1.** *Suppose that $(M, g)$ is a compact real analytic Riemannian manifold, with $L^2$-normalized eigenfunctions $\{\phi_\lambda\}$ satisfying $||\phi_\lambda||_\infty = \Omega(\lambda^{n-1})$. Then $(M, g)$ is a $Y_\ell^m$-manifold for some $\ell, m$. In particular, $M = S^2$ (topologically) if $n = 2$.*

We recall (cf. [Besse], Definition 7.8 (p. 182)) that a compact manifold $M$ with a distinguished point $m$ is called a $Y_\ell^m$-manifold if all geodesics issuing from $m$ come back to $m$ at time $\ell$. We further recall that a compact $Y_\ell^m$-manifold has the properties that $\pi_1(M)$ is finite and $H^*(M, \mathbb{Q})$ has exactly one generator.

**Proof of Theorem 5.1** By assumption, there exists a constant $C > 0$ and a subsequence of $L^2$-normalized eigenfunctions $\{\phi_{\lambda_j}\}$ such that $||\phi_{\lambda_j}||_\infty \geq C\lambda_j^{\frac{(n-1)}{2}}$. This contradicts the last statement of Theorem (2.1), hence there exists a point $m$ such that $\mu_m(\mathcal{L}_m) > 0$.

Since $g$ is real analytic, $\exp : T_m^* M \to M$ is a real analytic map, hence $\Gamma_m$ is an analytic set. In any local coordinate patch $U \subset M$ containing $m$, $\Gamma_m$ is the zero set of a pair of real-valued real-analytic functions, i.e. has the form $(f_1(\xi), \ldots, f_n(\xi)) = (m_1, \ldots, m_n)$. The solutions are the same as for $(f_1(\xi) - m_1)^2 + \cdots + (f_n(\xi) - m_n)^2 = 0$, so $\Gamma_m$ is the zero set of a real analytic function.

It is well-known (see [BM, H, L, S]) that the zero set of a real analytic function is locally a finite union of embedded real analytic submanifolds $Y_i^{k_i}$ of dimensions $1 \leq k_i \leq n-1$. Thus, for each $\xi$ there exists a ball $B_\delta(\xi)$ such that

(43) $$\Gamma_m \cap B_\delta(\xi) = \cup_{i=1}^d Y_i^{k_i}.$$

We claim that for some $(\xi, \delta)$ there exists a component $Y_i^{n-1}$ of dimension $n-1$. If not, $\Gamma_m$ is of Haussdorf dimension $\leq n-2$. But then its radial projection $\rho : \Gamma_m \to S_m^* M$ would also have Haussdorf dimension $\leq n-2$. In fact, each ray through the origin in $T_m^* M$ intersects $\Gamma_m$ in at most countably many points. So the radial projection preserves the dimension. But this contradicts the fact that $\Lambda_m = \rho(\Gamma_m)$ has positive Lebesgue measure.

Now let $\Xi \subset \Gamma_m$ be an open embedded real analytic hypersurface of $T^* M$. Consider the rays $\{t\xi : 0 \leq t \leq 1\} \subset T_m^* M$ and the union

$$C = \cup_{\xi \in \Xi} \{t\xi : 0 \leq t \leq 1\}.$$

Thus, each ray in $C$ exponentiates to a geodesic loop which returns at $t = 1$. As proved in Proposition 4.2, it follows by the first variation formula that the length $|\xi|$ of each loop must be a constant independent of $\xi \in \Xi$.

We conclude that $|\xi| = \ell$ for some $\ell \in \mathbb{R}^+$ and $\xi \in Y$. But this equation is real analytic, hence must hold on all of $\Xi$; hence $\Xi \subset \ell S_m^* M$. Again, by real analyticity, $\ell S_m^* M \subset \Gamma_m$. This is the same as saying that $(M, g)$ is a $Y_\ell^m$-manifold.

□

In §7, we shall prove that the converse to Theorem 1.2 holds automatically in the analytic setting.

6. Generic metrics: Proof of Theorem 1.4

We now prove that there exists a residual subset of the space $\mathcal{G}$ of $C^\infty$ metrics with the Whitney $C^\infty$ topology for which $\|\phi\|_\lambda = o(\lambda^{\frac{n-1}{2}})$. In view of Theorem 1.1, it suffices to prove:

**Lemma 6.1.** *There exists a residual set $\mathcal{R} \subset \mathcal{G}$ of metrics such that $|\mathcal{L}_x^g| = 0$ for every $x \in M$ when $g \in \mathcal{R}$.*

*Proof.* Recall that $|\mathcal{L}_x^g| = 0$ if and only if $\int_{S^x M} 1/L_g^*(x,\xi) d\xi = 0$, if $L_g^*(x,\xi)$ is the loop-length function defined before.

To make use of this, choose a coordinate patch $\Omega \subset M$ with coordinates $y = \kappa(x)$ ranging over an open subset of $\mathbb{R}^n$. We then fix $K \subset \kappa(\Omega)$ be compact and let

$$F(g) = \sup_{y \in K} \int_{S^{n-1}} \frac{d\sigma}{L_g^*(y,\xi)},$$

using the induced coordinates $\{y,\xi\}$ for $T^*\Omega \subset T^*M$. Here also, $d\sigma$ is the standard surface measure on $S^{n-1}$, and we are abusing notation a bit by letting $L_g^*(y,\xi)$ denote the pushforward of $L_g^*$ using $\kappa$. It then suffices to show that the set of metrics for which

$$\mathcal{G}_N = \{g : F(g) < 1/N\}$$

are open and dense.

Density just follows from the fact that $F(g) = 0$ for any non-Zoll real analytic metric. Such metrics are dense in $\mathcal{G}$.

The main step in proving that these sets are also open is to show that the function $f(g,y) = \int_{S^{n-1}} d\sigma/L_g^*(y,\xi)$ is upper semicontinuous on $\mathcal{G} \times K$. This holds since $1/L_g^*(y,\xi)$ is a positive, (locally) bounded upper semicontinuous function on $\mathcal{G} \times \Omega \times S^{n-1}$, if we equip $\mathcal{G}$ with the $C^3$-topology. Therefore, if $(g_j, y_j) \to (g, y)$, $y_j \in K$, we have

$$\sup_j \int_{S^{n-1}} \frac{d\sigma}{L_{g_j}^*(y_j,\xi)} \leq \int_{S^{n-1}} \sup_j \frac{d\sigma}{L_{g_j}^*(y_j,\xi)}$$
$$\implies \limsup_j \int_{S^{n-1}} \frac{d\sigma}{L_{g_j}^*(y_j,\xi)} \leq \int_{S^{n-1}} \limsup_j \frac{d\sigma}{L_{g_j}^*(y_j,\xi)} \leq \int_{S^{n-1}} \frac{d\sigma}{L_g^*(y,\xi)},$$

using the dominated convergence theorem and upper semicontinuity.

Now let us prove that the sets $\mathcal{G}_N$ are open. Let $g \in \mathcal{G}_N$. By definition of $F$, $f(g,y) < \frac{1}{N}$ for each $y \in K$. Since $f$ is upper-semicontinuous, the set $\{f < \frac{1}{N}\}$ is open, so there exist $\delta(y)$ such that $B_{\delta(y)}(y) \times B_{\delta(y)}(g) \subset \{f < \frac{1}{N}\}$, if $B_\delta(y)$ and $B_\delta(g)$ denote the $\delta$ Euclidean and $C^3$ balls of $y$ and $g$, respectively. Here, $g$ is fixed so we do not indicate the dependence of the $\delta$'s on it. As $y$ varies over $K$, the balls $B_{\delta(y)}(y)$ give an open cover of $K$, and by compactness there exists a finite subcover $\{B_{\delta(y_j)}(y_j), \ j = 1, \ldots, N\}$. Let $\delta = \min_j \delta(y_j)$, so that $B_\delta(y_j) \times B_\delta(g) \subset \{f < \frac{1}{N}\}$ for $j = 1, \ldots, N$. If $g' \in B_\delta(g)$, then $f(y, g') < \frac{1}{N}$ for all $y \in K$. Hence $F(g') < \frac{1}{N}$ and so $\mathcal{G}_N$ is open. □

## 7. Converse results: Proofs of Theorems 1.5 and 1.6

The purpose of this section is to prove Theorems 1.5 - 1.6. We shall need to go into somewhat more detail about loopsets.

We first consider the simplest converse in which the loopsets are assumed to be clean.

### 7.1. Theorem 1.5 for manifolds with clean loopsets.

**Theorem 7.1.** *Suppose that $(M, g)$ is a compact Riemannian manifold. Then:*

**(i)** *If $\Gamma_x^T \neq \emptyset$ is clean, and $\{T\}$ is isolated in $Lsp_x$, then $R(\lambda, x) = o(\lambda^{n-1}) \implies |S\Gamma_x^T| = 0$.*

**(ii)** *If $(M, g)$ has clean loopsets. Then*

$$\tag{44} R(\lambda, x) = o(\lambda^{n-1}) \implies |\mathcal{L}_x| = 0.$$

We begin the proof with some considerations which are valid for all compact Riemannian manifolds.

For each $(T, x)$ we introduce the sequence of probability measures

$$\tag{45} \mu_{T,\lambda,x}(\theta) = \frac{1}{E_\lambda(x,x)} \sum_{\lambda_\nu \leq \lambda} |\phi_j(x)|^2 \delta_{e^{iT\lambda_j}}(\theta)$$

on $S^1$. Here, $\delta_{e^{iT\lambda_j}}(\theta)$ denotes the point mass at the indicated point of $S^1$.

**Lemma 7.2.** *Let $(M, g)$ be any compact Riemannian manifold. Suppose that $R(\lambda, x) = o(\lambda^{n-1})$. Then, $\mu_{T,x,\lambda} \to d\theta$ for all $T$.*

*Proof.* The $k$th Fourier coefficient of $\mu_{T,\lambda,x}(\theta)$ equals

$$\tag{46} \hat{\mu}_{T,\lambda,x}(\theta) = \frac{1}{E_\lambda(x,x)} \sum_{\lambda_\nu \leq \lambda} |\phi_j(x)|^2 e^{iTk\lambda_j}.$$

Hence $R(\lambda, x) = o(\lambda^{n-1})$ implies

$$\tag{47} \begin{aligned} E_\lambda(x,x)\hat{\mu}_{T,\lambda,x}(\theta) &= \int_0^\lambda e^{ikT\lambda} dE_\lambda(x,x) \\ &= \int_0^\lambda e^{ikT\lambda} d\lambda^n + \int_0^\lambda e^{ikT\lambda} dR(\lambda, x) \\ &= n \int_0^\lambda e^{ikT\lambda} \lambda^{n-1} d\lambda + e^{ikT\lambda} R(\lambda, x) - ikT \int_0^\lambda e^{ikT\lambda} R(\lambda, x) \\ &= o(\lambda^n), \text{ for } k \neq 0 \end{aligned}$$

Since $d\theta$ is the unique probability measure $\mu$ on $S^1$ satisfying $\hat{\mu}(k) = \delta_{k0}$, we have

$$\tag{48} \int_{S^1} e^{ik\theta} d\mu_{T,\lambda,x}(\theta) \to 0 \quad (\forall k \neq 0) \iff \{d\mu_{T,\lambda,x}(\theta) \to d\theta\}.$$

Hence Lemma 7.2 follows from (47)-(48). □

In view of Lemma 7.2, to complete the proof it suffices to show:

**Lemma 7.3.** *We have:*

**(i)** *Suppose that* $\{T\}$ *is an isolated point in* $Lsp_x$ *and that* $\Gamma_x^T$ *is clean. Then* $\mu_{T,\lambda,x}(\theta) \to d\theta$ *implies that* $S\Gamma_x^T$ *is a submanifold of* $S_x^*M$ *of dimension* $\leq n-2$.

**(ii)** $\mu_{T,\lambda,x}(\theta) \to d\theta \;\; (\forall T \neq 0) \implies |\mathcal{L}_x| = 0$.

*Proof.* **(i)** By (48), (i) asserts that if $\Gamma_x^T$ is clean, then

$$(49) \qquad \int_{S^1} e^{ik\theta} d\mu_{T,\lambda,x}(\theta) \to 0 \;\; (\forall k) \implies |S\Gamma_x^T| = 0.$$

We prove this by studying the asymptotics of $[U(T)E_\lambda](x,x)$ as $\lambda \to \infty$, which are dual to the singularity of $U(T+s)(x,x)$ at $s=0$. These asymptotics were studied (in greater generality) in [Z] (§2), and we follow that discussion in the calculation to follow.

Following [Z], we write $U(T+s,x,x) = \int_{M \times M} U(T+s,y,y')\delta_x(y)\delta_x(y')dydy'$. Viewing $x$ and $T$ as fixed, we obtain a distribution in $s$ with singularities at times in the loop-length spectrum $Lsp_x$ (referred to as 'sojourn times' in the more general context of [Z]). If $\{T\}$ is isolated in $Lsp_x$, then there exists an interval around $s=0$ so that $U(T+s,x,x)$ is singular only at $s=0$. If $\Gamma_x^T$ is a clean intersection, it then follows from the clean composition calculus that $U(T+s,x,x)$ is a Lagrangean distribution in $s$ in this interval (see [Z], §2, for the specific case at hand, or [DG] [Ho IV] for the general theory). This implies that $U(T+s,x,x)$ has an isolated singularity at $s=0$ with a singularity expansion,

$$U(T+s,x,x) = a_r(T,x)(s+i0)^{-r} + a_{r-1}(T,x)(s+i0)^{-r+1} + \cdots$$

where $\cdots$ denotes smoother terms, and where the order $r = \dim S\Gamma_x^T$. Equivalently, we may write

$$(50) \qquad U(T+s)(x,x) = \int_{-\infty}^{\infty} \frac{\partial A_T(x,\lambda)}{\partial \lambda} e^{-is\lambda} d\lambda,$$

where $A_T(x,\lambda)$ is a symbol of order $r-1$ in $\lambda$. By Lemma 2.5 or by ([Z], Lemma 3.1), it follows that

$$(51) \qquad [U(T)E_\lambda](x,x) - A_T(x,\lambda) = O(\lambda^{n-1}).$$

Therefore, the statement (i) in Lemma 7.3 is equivalent to

$$(52) \qquad \text{order of} \;\; A_T(x,\lambda) \leq n-1 \iff |\Gamma_x^T| = 0.$$

The principal symbol of $U(T+s,x,x)$ at $s=0$ (for fixed $x$) was calculated in ([Z], Proposition (1.10)-Lemma(3.6)-Theorem (3.15)). The formula depends on $\dim S\Gamma_x^T$. When $\dim S\Gamma_x^T = n-1$, then the order of $A_T$ equals $n$ and the symbol is given by

$$(53) \qquad \sigma_n(A_T(x,\cdot)) = |S\Gamma_x^T|.$$

Now suppose that $d\mu_{T,x,\lambda} \to d\theta$. Then the first moment of $d\mu_{T,x,\lambda}$ tends to zero. If $\dim S\Gamma_x^T = n-1$, it must equal $S_x^*M$. But then the order $n$ symbol of $A_T$ cannot equal zero. This contradiction completes the proof of (i). $\square$

**(ii)** This case follows from (i) and is easier since the hypothesis is much stronger. By Proposition 4.2, if $|\mathcal{L}_x| > 0$, then there exists $T > 0$ such that $\Gamma_x^T = TS_x^*M$. As just proved, this implies that the symbol of order $n$ of $A_T$ is positive. But $R(x,\lambda) = o(\lambda^{\frac{n-1}{2}})$ implies that it equals zero. $\square$

7.2. **Real analytic metrics.** We also have a converse theorem for real analytic metrics, which is a simple corollary of the above theorem.

**Corollary 7.4.** *Let $(M, g)$ be real analytic. Then $R(x, \lambda) = o(\lambda^{n-1}) \implies |\mathcal{L}_x| = 0$.*

*Proof.* The set of vectors $\{\xi \in S_x^* : \exp_x \xi = x\}$ is an analytic set. Hence, it has only countably many components. Since $L^*$ is constant on smooth components, $Lsp_x$ is countable for each $x$.

Suppose now that $|\mathcal{L}_x| > 0$ for some $x$. Then, (as discussed in the proof of Theorem 5, there must exist $T > 0$ such that $S\Gamma_x^T = TS_x^*M$. The time $\{T\}$ is isolated in $Lsp_x$. Moreover, $\Gamma_x^T$ is necessarily a clean intersection, since by homogeneity of $\exp TH_p$ we have $(\exp TH_p)T_x^*M = T_x^*M$. Hence, the estimate follows from Theorem 7.1. □

7.3. **Almost Clean.** We now generalize the result to 'almost-clean' loopsets. In the argument above, we assumed that for $(T, x)$, $U(T+s)(x, x)$ is a Lagrangean distribution in $s$, but much less is sufficient since we only used the principal symbol and we only required $o(\lambda^n)$ from the remainder estimate.

**Definition 7.5.** *We say that $\Gamma_x^T$ is an 'almost-clean' loopset if, $\Gamma_x^T \subset T_x^*M$ is a (homogeneous) submanifold with boundary and if the intersection $\exp TH_p T_x^*M \cap T_x^*M$ is a clean intersection in the interior of $\Gamma_x^T$. We say that $(M, g)$ has almost clean loopsets if this property holds for all $T$.*

Simple examples of metrics with an almost clean loop-set $\Gamma_x^T$ are provided by generic bumps on standard spheres. In other words, let $B_x(r)$ denote the ball of radius $r$ around a point $x \in (S^n, g_0)$ (the standard sphere), and let $e^{2u}g_0$ denote a conformal perturbation with supp $u = B_x(r)$. For points $x$ in the exterior of the 'bump' $B_x(r)$, $\Gamma_x^{2\pi}$ contains any direction which avoids the bump. For generic $u$ it will not contain any other loops at time $T = 2\pi$. Hence, $S\Gamma_x^{2\pi}$ will be a submanifold with boundary of $S_x^*M$, with boundary the set of initial tangent vectors to geodesics which intersect $\partial B_x(r)$ tangentially. More generally, we could put several bumps on the sphere and obtain a number of components. Instead of bumps, one could put handles, and thus obtain any topology. This example will be important in the proof of Theorem 1.6.

Since $\Gamma_x^T$ could have several components, there is a possibility that cancellation between components can cause slow growth in $U(T)E_\lambda(x, x)$ even when $|\mathcal{L}_x| > 0$. We adopt a reasonable hypothesis, easy to check in our main examples, which rules this out.

**THEOREM 7.6. (i)** *Suppose that $\{T\}$ is an isolated point in $Lsp_x$ of $(M, g)$, that $\Gamma_x^T$ is an almost-clean loopset, and that the (common) Morse indices of loops in all components of $\Gamma_x^T$ of dimension $n-1$ are the same. Then $R(\lambda, x) = o(\lambda^{n-1}) \implies |S\Gamma_x^T| = 0$.*

**(ii)** *Suppose that $(M, g)$ has almost-clean loopsets, and that, for all $T$, the (common) Morse indices of loops in all components of $\Gamma_x^T$ of dimension $n-1$ are the same. Then*

$$R(\lambda, x) = o(\lambda^{n-1}) \implies |\mathcal{L}_x| = 0. \tag{54}$$

*Proof.* We only prove (i) since, as in the clean case, (ii) is an immediate consequence.

Nothing changes from the proof in the clean case until (50). Since the clean composition hypothesis is not assumed to hold, $U(T+s, x, x)$ is not known to be Lagrangean, although the singularity at $s = 0$ is still assumed to be isolated. We therefore cut the wave kernel up into pieces which can be controlled easily. The following argument is somewhat similar to the proof of [Z] Theorem (3.20).

We fix $(T, x)$ and choose a cover $S^*M$ adapted to $S\Gamma_x^T$ as follows. By assumption, $S\Gamma_x^T$ is a submanifold with boundary of $S^*M$ for all $T$. It therefore has only finitely many connected components. We write

$$S\Gamma_x^T = \cup_{j=1}^{n_T} K_x^{Tj} \cup \cup_{n=1}^{r_T} N_x^{Tj},$$

where $\{K_x^{Tj}\}$ are the components (possibly empty) of dimension $n-1$ and $\{N_x^{Tj}\}$ are the components of dimensions $\leq n-2$. The Morse indices of the loops in each component $K_x^{Tj}$ are the same; we are assuming the Morse index is the same constant for all $K_x^{Tj}$, as is the case in many natural examples.

We then form the cover $\{intK_x^{Tj}, T_\varepsilon(\partial K_x^{Tj}), T_\varepsilon(N_x^{Tj}), \mathcal{O}\}$, where $T_\varepsilon(\partial K_x^{Tj})$ is the tube (tubular neighborhood) of radius $\varepsilon$ around the boundary of $K_x^{Tj}$, where $T_\varepsilon(N_x^{Tj})$ is the $\varepsilon$-tube around $N_x^{Tj}$, and where $\mathcal{O}$ covers the rest of $S^*M$.

We then define a partition of unity subordinate to the cover as follows: we define functions $b_\varepsilon^j \in C_0^\infty(intK_x^{Tj})$ which are invariant under $\exp TH_p$ and $\equiv 1$ in points of $intK_x^{Tj}$ which lie outside of $T_\varepsilon(\partial K_x^{Tj})$. We further define $c_\varepsilon^j \in C_0^\infty(T_\varepsilon(\partial K_x^{Tj}))$, $d_\varepsilon^j \in C_0^\infty(T_\varepsilon(N_x^{Tj}))$, and $r_\varepsilon \in C_0^\infty(\mathcal{O})$ so that $\sum_{j=1}^{n_T} b_\varepsilon^j + \sum_{j=1}^{n_T} c_\varepsilon^j + \sum_{j=1}^{r_T} d_\varepsilon^j + r_\varepsilon \equiv 1$. We note that $b, c, d, r$ also depend on $T$; since it is fixed throughout the calculation, we do not indicate this in the notation.

In a well-known way (see [Ho IV]), we may quantize the symbols to obtain a pseudo-differential partition of unity:

$$\sum_{j=1}^{n_T} B_\varepsilon^j + \sum_{j=1}^{n_T} C_\varepsilon^j + \sum_{j=1}^{r_T} D_\varepsilon^j + R_\varepsilon \sim I.$$

We then have

$$U(T)E_\lambda(x,x) = [\sum_{j=1}^{n_T} B_\varepsilon^j U(T)E_\lambda](x,x) + [\sum_{j=1}^{n_T} C_\varepsilon^j U(T)E_\lambda](x,x) + [\sum_{j=1}^{r_T} D_\varepsilon^j U(T)E_\lambda](x,x)$$
$$+ [R_\varepsilon U(T)E_\lambda](x,x).$$

Since $B_\varepsilon^j$ is microsupported in the set where $\Gamma_x^T$ is a clean intersection, it follows that $B_\varepsilon^j U(T+s)](x,x)$ is a Lagrangean distribution for $s$ near $0$ and its principal symbol of order $n$ is given by $(i^{m_j} \int_{S\Gamma_x^T} b_\varepsilon^j d\nu_x)$ as in the clean case. Here, $m_j$ is the Morse index of loops in $K_x^{Tj}$ (which we assume independent of $j$). It follows by the standard trace asymptotics (see Lemma 2.5) that

$$[B_\varepsilon^j U(T)E_\lambda](x,x) = i^{m_j}(\int_{S\Gamma_x^T} b_\varepsilon^j d\nu_x)\lambda^n + O(\lambda^{n-1}).$$

Also, $R_\varepsilon U(T+s, x, x)$ is smooth is smooth in $s \in \mathbb{R}$ since no loops at $x$ occur in the microsupport of $R_\varepsilon$. We conclude that

$$(55) \quad U(T)E_\lambda(x,x) = \left(\sum_{j=1}^{n_T} \int_{S\Gamma_x^T} b_\varepsilon^j d\nu_x\right)\lambda^n + \left([\sum_{j=1}^{n_T} C_\varepsilon^j U(T)E_\lambda](x,x)\right) \\ + [\sum_{j=1}^{r_T} D_\varepsilon^j U(T)E_\lambda](x,x)) + O(\lambda^{n-1}).$$

We estimate the two middle terms in the same way; the harder one is the $C$-term, so we leave the details of the $D$ - terms to the reader. We have:

$$|C_\varepsilon^j U(T)E_\lambda](x,x)| = |\sum_{\lambda_\nu \leq \lambda} e^{i\lambda_\nu t} \phi_\nu(x) C_\varepsilon^j \phi_\nu(x)|$$

$$\leq \sum_{\lambda_\nu \leq \lambda} |\phi_\nu(x)||C_\varepsilon^j \phi_\nu(x)| \leq \sqrt{\sum_{\lambda_\nu \leq \lambda} |\phi_\nu(x)|^2} \sqrt{\sum_{\lambda_\nu \leq \lambda} |C_\varepsilon^j \phi_\nu(x)|^2}$$

(56)

$$\sim \lambda^n \sqrt{\int_{S_x^* M} |c_\varepsilon^j|^2 d\nu_x}$$

$$= O(\sqrt{Vol(T_\varepsilon(\partial K_x^{Tj}))}\lambda^n) = O(\sqrt{\varepsilon}\lambda^n).$$

In the above we used that, for any 0th order pseudodifferential operator,

$$\sum_{\lambda_\nu \leq \lambda} |C_\varepsilon \phi_\nu(x)|^2 = CE_\lambda C^*(x,x) = \left(\int_{S_x^* M} |c|^2 d\nu_x\right)\lambda^n + O(\lambda^{n-1}).$$

Since $T_\varepsilon(\partial K_x^{Tj})$ is a tube around a hypersurface in $S^*M$, the integral is $O(\varepsilon)$. Thus, by the same Fourier coefficient calculation as in the clean case, $R(\lambda, x) = o(\lambda^{n-1})$ implies that $U(T)E_\lambda(x,x) = o(\lambda^n)$ for all $T \neq 0$. It follows for any $\varepsilon > 0$ that

$$\sum_{j=1}^{n_T} i^{m_j} \int_{S\Gamma_x^T} b_\varepsilon^j d\nu_x = O(\sqrt{\varepsilon}).$$

By the assumption that $m_j \equiv m$ for all $j$, it follows then that $|S\Gamma_x^T| = 0$. This completes the proof of (i).

To prove (ii), it suffices to observe that $Lsp_x$ is countable, and we have just seen that each $S\Gamma_x^T$ has measure zero. □

### 7.4. Counterexamples: Proof of Theorem 1.6.
We now present an example of a $C^\infty$ (but not analytic) surface $(M, g)$ satisfying:

$$\exists x : |\mathcal{L}_x| > 0, \ |R(\lambda, x)| = \Omega(\lambda^{\frac{n-1}{2}}) \text{ but } L^\infty(\lambda, g) = o(\lambda^{\frac{n-1}{2}}).$$

Thus, neither $|\mathcal{L}_x| > 0$ nor $|R(\lambda, x)| = \Omega(\lambda^{\frac{n-1}{2}})$ is sufficient to imply maximal eigenfunction growth.

Our example is of a surface of revolution, i.e. a surface $(M, g)$ whose metric is invariant under an action of $S^1$ by isometries. The two possible $M$ are $M = S^2$ and $M = T^2$ (the two-torus). On $M = S^2$ we can construct metrics such that $\exists x : |\mathcal{L}_x| > 0$, but no sequence of eigenfunctions has maximal growth at $x$; but the invariant eigenfunctions will have maximal growth at the poles. To obtain an example of a surface which does not have eigenfunctions of maximal growth, we instead choose $M = T^2$.

We often write $M = S^1 \times S^1$ with the angle coordinates $(x, \theta)$. We shall refer the circle over which $x$ varies as the base circle, and the orbits as the fiber circle. Our torus of revolution is of the following kind:

**Definition 7.7.** $(T^2, g)$ *will be called a torus of revolution with an equatorial band if in the standard coordinates $(x, y)$ on $T^2$ we have $g_a = dx^2 + a(x)^2 d\theta^2$, where $a = \sqrt{1 - x^2}$ for $x \in (-\varepsilon, \varepsilon)$.*

Such a surface may be obtained by by revolving the graph $y = e^{u(x)}$, $0 \leq x \leq 2\pi$ of a $2\pi$-periodic strictly positive function $e^{u(x)}$ around the $x$-axis, and then joining the ends $x = -1, x = 1$ together. The 'equatorial band' is the part $B_\varepsilon \sim (\varepsilon, \varepsilon) \times S^1$ lying over angular interval $(-\varepsilon, \varepsilon)$ (for some $0 < \varepsilon < 1/2$), and is isometric to an equatorial band lying over the same interval of the $z$-axis of the standard $S^2$. Clearly, $|\mathcal{L}_x| = 2\pi$ for the open annulus of $x \in B_\varepsilon$.

We now construct $u$ such that the eigenfunctions of $(T^2, g_u)$ fail to have maximal growth. We break up the base circle $x \in S^1$ into four parts: a flat part on $[-1, -2\varepsilon] \cup [2\varepsilon, 1]$, a spherical part $[-\varepsilon, \varepsilon]$ and a 'bridge' over $(-2\varepsilon, -\varepsilon) \cup (\varepsilon, 2\varepsilon)$ between them. We shall construct a metric which is invariant under the involution $x \to -x$ (i.e. $e^{ix} \to e^{-ix}$), as follows: Introduce $C^\infty$ cutoff functions on $S^1$ satisfying:

$$\chi_1(x) = \begin{cases} 1, & x \in [-\varepsilon, \varepsilon] \\ 0, & x \in [-1, -2\varepsilon] \cup [2\varepsilon, 1] \end{cases}$$

$$\chi_2(x) = \begin{cases} 0, & x \in [-\varepsilon, \varepsilon] \\ 1, & x \in [-1, -2\varepsilon] \cup [2\varepsilon, 1] \end{cases}$$

We now define a class of profile curves of the form:

$$a_f(x) = \sqrt{1 - x^2} \chi_1(x) + \frac{1}{2}\sqrt{1 - \varepsilon^2} \chi_2(x) + f(x)(1 - \chi_1(x) - \chi_2(x)),$$

$$f(x) = f(-x), \quad f \in C^\infty([\varepsilon, 2\varepsilon]).$$

For concreteness, we shall construct $f$ so that $a_f$ is monotonically decreasing from the boundary of the round region (where $a_f(\varepsilon) = \sqrt{1 - \varepsilon^2}$) to the flat region (where $a_f \equiv \frac{1}{2}\sqrt{1 - \varepsilon^2}$). We shall begin with one such function $f_0$ and then consider small perturbations $f = f_0 + \varepsilon \dot{u}$. We denote the associated surface of revolution by $(T^2, g_f)$, and its Laplacian by $\Delta_f$. It is clear that only the values of $f$ on $\pm([\varepsilon, 2\varepsilon])$ will matter. We equip the space of $\dot{u} \in C^\infty([\varepsilon, 2\varepsilon])$ with the Whitney $C^\infty$-topology.

Our first result involves the splitting of eigenspaces $E_\lambda = \oplus_{n=1}^{mult(\lambda)} V_{\lambda, n}$ of $\Delta_{f_0}$ into irreducibles for $S^1$. Here, $V_{\lambda, n}$ denotes the joint eigenspace of $\Delta, \frac{\partial}{\partial \theta}$. Because the eigenvalues of Sturm-Liouville operators on $S^1$ have multiplicity at most two, and $\dim V_{\lambda, n} = 2$, the multiplicity of each $V_\lambda$ in $E_\lambda$ is at most one. Therefore, $mult(\lambda)$ counts the total number of distinct irreducibles in $E_\lambda$. The following Lemma shows that by a generic small perturbation of $f_0$, we can arrange that all of the eigenspaces are irreducible. We note that the same proof shows more generally that all eigenspaces are irreducible for generic surfaces of revolution. This result is analogous to the earlier result of one of the authors [Z3] that all eigenspaces are irreducible for generic $G$-invariant metrics on Riemannian $d$-manifolds when $G$ is a finite group and $d \geq d_G$ for some minimal dimension $d_G$ depending on $G$. We refer to [Z3] for background, references and arguments that are valid when $G = S^1$.

**Lemma 7.8.** *Let $\mathcal{M}_{irr}$ denote the subspace of $C^\infty([\varepsilon, 2\varepsilon])$ of functions $f$ such that all eigenspaces of $\Delta_f$ are irreducible for the $S^1$ action. Then $\mathcal{M}_{irr}$ is residual in $C^\infty([\varepsilon, 2\varepsilon])$. In particular, its intersection with any closed ball $B_\delta(f_0)$ around $f_0$ in $C^\infty([\varepsilon, 2\varepsilon])$ is residual in the ball.*

*Proof.* We need to prove that $\mathcal{M}_{irr}$ contains a countable intersection of open dense sets. As in [Z3], we denote by $\mathcal{M}_k$ the metrics for which all of the first $k$ eigenspaces are irreducible. It suffices to prove that each $\mathcal{M}_k$ is open-dense.

Since eigenvalues and eigenfunctions move continuously with perturbations, it is a simple fact based on the minimax characterization of eigenvalues that $\mathcal{M}_k$ is open in $C^\infty([\varepsilon, 2\varepsilon])$ and hence that its intersection with any closed ball is a relatively open subset of the ball. We refer to [Z3] for the standard argument. The novel point to prove is the density of each $\mathcal{M}_k$.

We now separate variables. For each $n \in \mathbb{Z}$, we denote the joint eigenvalues and normalized eigenfunctions of $\Delta, \frac{\partial}{\partial \theta}$ by $\{\lambda_{n,j}^2, \Phi_{n,j}\}$. Thus, $\Phi_{n,j}(x, \theta) = e^{in\theta} \Phi_{n,j}(x, 0)$. We shall refer to such joint eigenfunctions as 'equivariant eigenfunctions'. With no loss of generality, we may conjugate $\Delta$ by $\sqrt{a}$ to put it into Liouville normal form (i.e. transform it to the half-density Laplacian). The functions $\Phi_{n,j}(x, 0)$ are then complex eigenfunctions of the associated Sturm-Liouville operator

$$L_n = \frac{d^2}{dx^2} + \frac{n^2}{a(r)^2} + V, \quad V = a^{-1/2} \frac{d}{dx} a^{-1} \frac{d}{dx} a^{-1/2}.$$

When $n \neq 0$, we have

$$\Phi_{n,\lambda}(x, \theta) = e^{in\theta}(u_{n,\lambda}(x) + iv_{n,\lambda}(x)),$$

where $u_{n,\lambda}, v_{n,\lambda}$ are the real eigenfunctions of $L_n$. When $n = 0$, there is only a one-dimensional space of invariant eigenfunctions.

To prove that $\mathcal{M}_k$ is dense, it suffices to prove that for each of first $k$ eigenspaces, it is possible to find a perturbation which splits the eigenspace into irreducibles. Since only a finite number of irreducibles are involved, it suffices to do so for one eigenspace at a time and for one pair of irreducibles in its decomposition. The case where $n = 0$ for one of the eigenspaces is somewhat different, so we assume $n, m \neq 0$ at first.

We argue by contradiction. Suppose that, for a metric $g_0$, the eigenspace $E_\lambda$ contains at least two complex normalized eigenfunctions $\Phi_{n,\lambda}, \Psi_{m,\lambda}$ with $m, n \neq 0$. We would like to find a curve of metrics which splits up their eigenvalues. It is clear that $m \neq n$ since the eigenvalues of Sturm-Liouville operators are at most double. We recall that the variation of the eigenvalue under a deformation $e^{2u_\varepsilon} g_0$ is given by

$$\dot{\lambda} = \langle \dot{\Delta} \Phi_{n,\lambda}, \Phi_{n,\lambda} \rangle = \langle (\dot{u}\Delta + \Delta \dot{u}) \Phi_{n,\lambda}, \Phi_\lambda \rangle = 2\lambda \langle \dot{u} \Phi_{n,\lambda}, \Phi_{n,\lambda} \rangle,$$

where $\dot\Delta$ is the variation of $\Delta$. The same calculation holds for $\Psi_{m,\lambda}$. If the eigenvalue does not split under any deformation, then we have

$$\langle \dot u \Phi_{n,\lambda}, \Phi_{n,\lambda}\rangle = \langle \dot u \Psi_{m,\lambda}, \Psi_{m,\lambda}\rangle, \quad (\forall \dot u \in C^\infty(\varepsilon, 2\varepsilon))$$

Since

$$\int_0^T (|\Phi_{n,\lambda}(x,0)|^2 - |\Psi_{m,\lambda}(x,0)|^2)\dot u dx = 0, \quad (\forall \dot u \in C^\infty(\varepsilon, 2\varepsilon))$$

it follows that $|\Phi_{n,\lambda}(x,\theta)| = |\Psi_{m,\lambda}(x,\theta)|$ in $(\varepsilon, 2\varepsilon) \times S^1$. It follows that on this region $\Phi_{n,\lambda}(x,\theta) = e^{i(n-m)\theta}e^{i\alpha(x)}\Psi_{m,\lambda}(x,\theta)$ for some real-valued function $\alpha$.

To analyze the situation, we write

$$\Phi(x,\theta) = e^{in\theta}e^{\rho(x)+i\tau(x)}.$$

Then $\Delta e^{in\theta}e^{\rho(x)+i\tau(x)} = \lambda e^{in\theta}e^{\rho(x)+i\tau(x)}$ gives

$$\rho'' + i\tau'' + ((\rho')^2 + 2i\tau'\rho' - (\tau')^2) + \frac{a'}{a}(\rho' + i\tau') = \lambda + \frac{n^2}{a^2}.$$

Again separating into real and imaginary parts this becomes:

(57) $$\begin{cases} \Re: \quad \frac{n^2}{a(x)^2} - \lambda^2 = \rho'' + (\rho')^2 - (\tau')^2 + \frac{a'}{a}\rho' \\ \\ \Im: \quad 0 = 2\rho'\tau' + \tau''(x) + \frac{a'(x)}{a(x)}\tau'(x). \end{cases}$$

From the imaginary part equation we obtain

(58) $$\tau' = e^{\tau'(0)}(a^{-1}e^{-2\rho}) \implies \tau = \tau(0) + e^{\tau'(0)}\int_0^x (a^{-1}e^{-2\rho})dy.$$

Thus, the phase $\tau$ is determined up to two constants $\tau(0), \tau'(0)$ by the modulus $e^\rho$ of $\Phi_{n,\lambda}$.

We further observe by combining this with the real part equation that $\rho$ solves:

$$\rho'' + (\rho')^2 + \frac{a'}{a}\rho' - e^{2\tau'(0)}(a^{-2}e^{-4\rho}) = \frac{n^2}{a(x)^2} - \lambda^2.$$

Suppose now that $\rho_1 = |\Phi|, \rho_2 = |\Psi|$. Subtracting their two equations gives

$$\frac{d^2}{dx^2}(\rho_1 - \rho_2) = (\rho_1' + \rho_2')\frac{d}{dx}(\rho_1 - \rho_2) + \frac{a'}{a}\frac{d}{dx}(\rho_1 - \rho_2) - a^{-2}e^{2\tau(0)}(e^{-4\rho_1} - e^{-4\rho_2}).$$

Here, we assumed that $\tau_1(0) = \tau_2(0)$, which is no loss of generality since we may multiply $\Psi$ by a $e^{i(\tau_1(0)-\tau_2(0))}$ obtain a scalar multiple with this property.

If we put $f = \rho_1 - \rho_2$, then we obtain the inequality

$$\tfrac{d^2}{dx^2}f = (\rho_1' + \rho_2' + \tfrac{a'}{a})f' - a^{-2}e^{2\tau(0)}fG(\rho_1, \rho_2)$$

$$\implies |f''(x)| \leq A|f'(x)| + B|f(x)|$$

for some constants $A, B$, namely the supremums of $(\rho_1' + \rho_2' + \frac{a'}{a})$ and $G$ respectively. It follows then by Aronsajn's unique continuation theorem that $f = 0$ on an interval implies $f \equiv 0$. Hence we have that $|\Phi(x)| = |\Psi(x)|$ globally.

Of course, this is no contradiction: there do exist metrics on $T^2$ possessing pairs equivariant eigenfunctions with the same eigenvalue and same pointwise norm: indeed, for the standard $\mathbb{R}^2/\mathbb{Z}^2$ one may take $\Phi_{n,\lambda}(x,\theta) = e^{i(kx+n\theta)}$, $\Psi_{m,\lambda}(x,\theta) = e^{i(\ell x+m\theta)}$, where $k^2 + n^2 = m^2 + \ell^2$. We now assert that such pairs only exist for flat metrics of revolution on $T^2$.

Indeed, suppose that $\Phi = F\Psi$. We also note that $F(x,\theta) = e^{i(n-m)\theta}F(x,0)$. Write $F(x,0) = f(x), \Psi(x,0) = \psi(x)$. Since also $|f(x)| = 1$, there exist smooth functions $\rho, \tau, \alpha$ such that
$$\Psi(x,\theta) = e^{im\theta}\psi(x), \quad \psi(x) = e^{\rho(x)+i\tau(x)}, \quad f(x) = e^{i\alpha(x)},$$
$$-\lambda^2 F\Psi = \Delta(F\Psi) = \Psi\Delta F + 2\nabla\Psi \cdot \nabla F + F\Delta\Psi \iff \Psi\Delta F + 2\nabla\Psi \cdot \nabla F = 0$$
$$\iff \nabla(\Psi^2 \nabla F) = 0.$$

Simple calculations give:
$$\nabla(\Psi^2 \nabla F) = 0 \implies \frac{m^2-n^2}{a(x)^2}f(x) = 2\frac{(\psi(x)^2)'}{\psi(x)^2}f'(x) + (f''(x) + \frac{a'}{a}f')$$
$$\implies \frac{m^2-n^2}{a(x)^2} = 2(\rho'(x) + i\tau'(x))[i\alpha'(x)] - \alpha'(x)^2 + i\alpha''(x) + i\frac{a'(x)}{a(x)}\alpha'(x).$$

Breaking up into real and imaginary parts gives
$$\begin{cases} \Re: & \frac{m^2-n^2}{a(x)^2} = -2\tau'(x)\alpha'(x) - (\alpha'(x))^2 \\ \\ \Im: & 0 = 2\rho'(x)\alpha'(x) + \alpha''(x) - \frac{a'(x)}{a(x)}\alpha'(x). \end{cases}$$

The imaginary part equation has the unique solution
$$\alpha'(x) = (n-m)a^{-1}(x)e^{-2\rho(x)}$$
consistent with the given initial conditions. Substituting into the real equation gives:
$$(59) \qquad m^2 - n^2 = -2(n-m)\tau'(x)a^{-1}(x)e^{-2\rho(x)} - (n-m)^2 e^{-4\rho(x)}.$$

Substituting (58) into (59) gives
$$(60) \qquad m^2 - n^2 = Ce^{-4\rho} \implies \rho(x) = C_{m,n}$$
for some constant $C_{m,n}$. Substituting this into the real part equation of (57) and substituting the value of $\tau'$ in (58) gives
$$(61) \qquad (C^2 - n^2)a^{-2} = \lambda,$$
which implies that $a(x)$ is constant. $\square$

To complete the proof we use the following Lemma.

**Lemma 7.9.** *For any 'bridge function' $f$, the joint $L^2$-normalized eigenfunctions $\{\Phi_{n,\lambda}\}$ of $\Delta$ and of the $S^1$ action on $(T^2, g_a)$ satisfy $\|\Phi_{n,\lambda}\|_\infty = o(\lambda^{\frac{n-1}{2}})$.*

*Proof.* Let us now prove eigenfunction estimates on our torus of revolution $T$ with metric
$$(62) \qquad g_a = dx^2 + (a(x))^2 d\theta^2, \quad -1 \leq x \leq 1, \; -\pi \leq \theta \leq \pi,$$
where of course $a$ is a positive periodic function on $\mathbb{R}/\mathbb{Z}$. We shall show that for $L^2$-normalized eigenfunctions of the form
$$(63) \qquad \Phi = \Phi_{n,\lambda} = e^{in\theta}\Phi(x)$$

with eigenvalue $\lambda$ we have the bounds

$$\|\Phi\|_{L^\infty(T^2)} \leq C\lambda^{1/2-\sigma}, \tag{64}$$

for some $\sigma > 0$ to be specified later. These bounds of course imply that for our the above metric we have $L^\infty(\lambda, g) = o(\lambda^{1/2})$; however, as we shall sketch after the elementary proof, the bounds we obtain from this simple argument are far from optimal.

To prove (64), fix an even function $\rho \in C_0^\infty(\mathbb{R})$ with integral one and having the property that

$$\rho(t) = 0 \quad \text{if } |t| > 1.$$

Then, if $\hat\rho$ is the Fourier transform of $\rho$, it follows that for a given $\delta > 0$

$$\Phi = \chi_\lambda^\delta \Phi = \hat\rho(\delta(\lambda - \sqrt{-\Delta}))\Phi. \tag{65}$$

Note that

$$\chi_\lambda^\delta = \frac{1}{2\pi\delta} \int e^{-it(\lambda - \sqrt{-\Delta})} \rho(\delta^{-1}t)\,dt. \tag{66}$$

Let $\chi_\lambda^\delta(w, z)$, $w, z \in T^2$, be the kernel of $\chi_\lambda^\delta$. Then by finite propagation speed for $\cos(t\sqrt{-\Delta})$ it follows that

$$\chi_\lambda^\delta(w, z) = 0 \quad \text{if } \operatorname{dist}(w, z) > 2\delta,$$

with $\operatorname{dist}(w, z)$ being the distance between $w$ and $z$ measured by the metric $g_a$. Because of this, we can explain the role of $\delta$ in our estimates. First of all, if $\delta > 0$ is small the kernel is supported on small scales, and so after rescaling to the unit scale in local coordinates behaves increasingly like ones arising from a flat metric as $\delta \to 0$. This makes calculations in the proof easier; however, as $\delta \to 0$ one sees from (65) that $\chi_\lambda^\delta$ looks less and less like an eigenspace projection operator, which will reflect negatively on (64). Ideally, one would like to take $\delta$ to be, say, half the injectivity radius, but then the calculations that arise in the proof seem formidable.

To proceed, let us note that if $\delta$ is larger than $\lambda^{1-\alpha_0}$, where $\alpha_0 > 0$ is fixed and if $\delta$ is smaller than half the injectivity radius then we can write

$$\chi_\lambda(w, z) = \lambda^{1/2}\delta^{-3/2} \sum_\pm a_{\lambda,\delta}^\pm(w, z) e^{\pm i\lambda \operatorname{dist}(w, z)},$$

where the amplitudes satisfy the bounds

$$|D_{w,z}^\gamma a_{\lambda,\delta}^\pm(w, z)| \leq C_\gamma \delta^{-|\gamma|},$$

with $C_\gamma$ independent of $\delta$. This follows from routine stationary phase calculations using (66). (The case $\delta \approx 1$ was carried out in [So3], p. 138-140.) Also, the support properties of $\chi_\lambda^\delta(w, z)$ imply that

$$a_{\lambda,\delta}^\pm(w, z) = 0 \quad \text{if } \operatorname{dist}(w, z) > 2\delta.$$

To proceed, let us fix $w \in T^2$ and use the natural coordinates $(x, \theta)$ corresponding to (62). In proving bounds for $\phi$ at $w$, in view of (63), we may assume that $w$ has coordinates $(x_0, 0)$. If $(x, \theta)$ then are the coordinates of $z$, it follows from (63) that

$$\phi(x_0, 0) = \sum_\pm \int K_{\lambda,n,\delta}^\pm(x_0, x)\phi(x)\,dx, \tag{67}$$

where $K^\pm_{\lambda,n,\delta}$ is given by the oscillatory integral

$$K^\pm_{\lambda,n,\delta}(x_0, x) = \lambda^{1/2} \delta^{-3/2} \int e^{\pm i\lambda \text{dist}((x_0,0),(x,\theta))} e^{in\theta} a^\pm_{\lambda,\delta}(x_0, 0, x, \theta) \, d\theta.$$

By applying the Schwarz inequality we of course get

$$|\Phi(x_0, 0)| \leq \sum_\pm \Bigl(\int |K^\pm_{\lambda,n,\delta}(x_0, x)|^2 \, dx\Bigr)^{1/2} \Bigl(\int |\Phi(x)|^2 \, dx\Bigr)^{1/2}$$

$$\leq C \sum_\pm \Bigl(\int |K^\pm_{\lambda,n,\delta}(x_0, x)|^2 \, dx\Bigr)^{1/2},$$

using in the last step the fact that $d\text{Vol} = a(x) dx d\theta$ with $a(x) > 0$ and the fact that $\Phi(x, \theta)$ has $L^2$-norm one.

Because of this our proof of (64) boils down to estimating the $L^2$ norms of $K^\pm_{\lambda,n,\delta}$ and then optimizing in $\delta$. To do this note that the integrand of the oscillatory integral defining $K^\pm_{\lambda,n,\delta}$ vanishes unless $|\theta| + |x_0 - x| \leq C\delta$. Therefore, it is natural to make the change of variables $\theta \to \delta\theta$, $x \to x_0 + \delta y$, which results in

(68) $$K^\pm_{\lambda,n,\delta}(x_0, x_0 + \delta y) = (\lambda/\delta)^{1/2} \int e^{\pm i\lambda\delta \Psi_\delta(x_0; y, \theta)} a^\pm_{\lambda,\delta}(x_0, 0; x_0 + \delta y, \delta\theta) \, d\theta,$$

where

$$\Psi_\delta = \delta^{-1} \text{dist}\bigl((x_0, 0), (x_0 + \delta y, \delta\theta)\bigr) + n\delta\theta.$$

If the warping factor $a(x)$ for the metric (62) were constant in a $2\delta$ neighborhood of $x_0$ then one would have $(\partial/\partial\theta)^2 \Psi_\delta = -(a(x_0)y)^2/|y^2 + (a(x_0)\theta)^2|^{3/2}$, which results in a lower bound of $(a(x_0)y)^2$ for the Hessian of the phase function in this case. For non-constant $a$ one has

$$\delta^{-1} \text{dist}\bigl((x_0, 0), (x_0 + \delta y, \delta\theta)\bigr) = \sqrt{y^2 + (a(x_0)\theta)^2 + r_\delta(x_0, y, \theta)},$$

on the support of the oscillatory integral where for where for a given $j = 0, 1, 2, \ldots$

$$|(\partial/\partial\theta)^j r_\delta| \leq C\delta.$$

Consequently, there are uniform constants $c_0 > 0$ and $C_0 < \infty$ so that

$$|(\partial/\partial\theta)^2 \Psi_\delta| \geq c_0 y^2, \text{ if } y^2 \geq C_0 \delta.$$

Since the amplitude in (68) has uniformly bounded derivatives, we conclude from stationary phase that

$$|K^\pm_{\lambda,n,\delta}(x_0, x_0 + \delta y)| \leq C(\lambda/\delta)^{1/2}(y^2 \delta \lambda)^{-1/2} = C\delta^{-1}|y| \text{ if } |y| \geq C_0 \delta^{-1/2}.$$

One also has the trivial bound $|K^\pm_{\lambda,n,\delta}| \leq C(\lambda/\delta)^{1/2}$. Using these two bounds we can compute the $L^2$ norm of the kernel $K^\pm_{\lambda,n}$ with respect to the original coordinates:

$$\Bigl(\int |K^\pm_{\lambda,n,\delta}(x_0, x)|^2 \, dx\Bigr)^{1/2} = \delta^{1/2} \Bigl(\int |K^\pm_{\lambda,n,\delta}(x_0, x_0 + \delta y)|^2 dy\Bigr)^{1/2}$$

$$\leq C\delta^{-1/2} \Bigl(\int_{\delta^{1/2}}^{+\infty} y^{-2} \, dy\Bigr)^{1/2} + C\lambda^{1/2} \Bigl(\int_0^{\delta^{1/2}} dy\Bigr)^{1/2}$$

$$\leq C\delta^{-3/4} + C\lambda^{1/2}\delta^{1/4}.$$

One optimizes this by choosing $\delta = \lambda^{-1/2}$, which implies that one has (64) with $\sigma = 1/8$. □

**Completion of the proof of Theorem 1.6**

Two of the three properties follows from the corollary and from the fact that $|\mathcal{L}_x| > 0$ within the band. For the third property, we need to prove the existence of a point $x_0$ such that $R(\lambda, x_0) = \Omega(\lambda^{\frac{n-1}{2}})$. For this we appeal to Theorem 7.6 (i). Let $x$ be a point on the equator of the round part. Then $\Gamma_x^{2\pi}$ has an almost-clean contribution coming from the great circles in the equatorial band. We claim that by choosing the flat part to be long enough, there $\Gamma_x^{2\pi}$ contains no other loops. Indeed, any other loop would have to enter the flat part. It cannot reverse direction, so it must travel around the flat part before re-entering the spherical part. As long as the flat part is at least $2\pi$ units in length, any such geodesic loop at $x$ will not belong to $\Gamma_x^{2\pi}$. □

We believe that the correct bound should be $\lambda^{1/4}$. We shall discuss this further below.

## 8. Further problems and conjectures

Let us mention some related natural problems which remain open.

**Problem 1** First, what is a sufficient condition for maximal eigenfunction growth at a given point $x$ or at some point of $(M, g)$?

As mentioned above, $|\mathcal{L}_x| > 0$ is not sufficient to imply the existence of a sequence of eigenfunctions blowing up at $x$ at the maximal rate. In fact, to our knowledge, the only Riemannian manifolds known to exhibit maximal eigenfunction growth are surfaces of revolution and compact rank one symmetric spaces. In the case of surfaces of revolution, invariant eigenfunctions must blow up at the poles because all other eigenfunctions vanish there and yet the local Weyl law (4) must hold. In the case of compact rank one symmetric spaces, the exceptionally high multiplicity of eigenspaces allows for the construction of eigenfunctions of maximal sup norm growth. In each case, $\mathcal{L}_x$ has full measure, but the mechanisms producing maximal eigenfunction growth involve something more. All of these examples are completely integrable and eigenfunctions with maximal sup-norms may be explicitly constructed by the WKB method. Let us describe the symplectic geometry, because it turns out that we are tantalizingly close to proving it must occur in the real analytic case.

The eigenfunctions with maximal sup norms in these cases are actually oscillatory integrals (quasimodes) associated to (geodesic flow)-invariant Lagrangean submanifolds diffeomorphic to $S_x^* M \times S^1 \sim S^{n-1} \times S^1$. They are the images of $S^1 \times S_m^* M$ under the Lagrange immersion

$$\iota : S^1 \times S_m^* M \to T^* M \backslash 0, \quad \iota(t, x, \xi) = G^t(x, \xi),$$

where $G^t : T^* M \backslash 0 \to T^* M \backslash 0$ is the geodesic flow. Under the natural projection $\pi : S^* M \to M$, $\pi : \iota(S^1 \times S_m^* M) \to M$ the sphere $S_x^* M$ is 'blown down' $S^* M$ to $m$. As discussed in [TZ] (and elsewhere), singularities of projections of Lagrangean submanifolds correspond closely to sup norms of the associated quasimodes, and such blow-downs give the maximal growth rate of the associated quasimodes. Existence of a quasimodes of high order attached to an invariant Lagrangean $S_x^* M \times S^1$ may therefore be the missing condition.

In the real analytic case, we show in Theorem 5.1 that all geodesics issuing from some point $m$ return to $m$ at a fixed time $\ell$. If, as is widely conjectured, they are smoothly closed curves at $m$, then these geodesics fill out a Lagrangean submanifold of the form $\iota(S^1 \times S_m^* M)$ above. To complete the conjectured picture, we would to show that there exist eigenfunctions which are oscillatory integrals associated to this Lagrangean. It is automatic that quasimodes (approximate eigenfunctions) can be constructed, but it is not necessarily the case that they approximate actual eigenfunctions.

We do not expect such quasimodes always to approximate eigenfunctions, nor do we expect maximal eigenfunction growth in all situations where $\mathcal{L}_x$ has full measure, even at all points $x$. In other words, we do not expect maximal eigenfunction growth on all Zoll manifolds (manifolds all of whose geodesics are closed). It is quite conceivable that (non-rotational) Zoll spheres do not have maximal eigenfunction growth, even though the converse estimate $R(\lambda, x) = \Omega(\lambda^{\frac{n-1}{2}})$ holds. Indeed, it is known ([V]) that on the standard $S^2$, almost every orthonormal basis of eigenfunctions satisfies $||\phi_\lambda|| = O(\sqrt{\log \lambda})$, even though special eigenfunctions (zonal spherical harmonics) have maximal eigenfunction growth. It is possible that eigenfunctions of typical Zoll surfaces resemble such typical bases of spherical harmonics rather than the special ones with extremal eigenfunction growth. In short, the converse direction appears to be a difficult open problem.

In a future article with J. Toth, we plan to study $L^p$-norms of eigenfunctions and quasimodes directly using oscillatory integral and WKB formulae. As a very special case, we will improve the sup norm bound on eigenfunctions of tori of revolution to $\lambda^{1/4}$.

At the opposite extreme is the problem of characterizing compact Riemannian manifolds with minimal eigenfunction growth, such as occurs on a flat torus. In [TZ] it is proved that in the integrable case, the only examples are flat tori and their quotients.

**Problem**: Characterize $(M, g)$ with maximal $L^p$-norms of eigenfunctions.

We have left this problem open for $p \leq \frac{2(n+1)}{n-1}$, and we expect the condition on $(M, g)$ to change at the the critical Lebesgue exponent $p = \frac{2(n+1)}{n-1}$. Indeed, as one of the authors proved in [So2], the 'geometry' of extremal eigenfunctions changes at this exponent. Specifically (cf. [So3], p. 142-144, [So1]), for $p > \frac{2(n+1)}{n-1}$ eigenfunctions concentrated near a point tend to have extreme $L^p$ norms, while for $2 < p < \frac{2(n+1)}{n-1}$ ones concentrated along stable closed geodesics tend to have this property; if $p = \frac{2(n+1)}{n-1}$, at least in the case of the round sphere, both types give rise to $\Omega(\lambda^{\delta(p)})$ bounds.

Bourgain [B] has constructed a metric on a 2-torus of revolution for which the maximal $L^6$ bound is attained although $|\mathcal{L}_x| = 0$ for all $x$. The eigenfunctions are similar to the highest weight spherical harmonics on $S^2$ which concentrate on the equator. Note that this is the Lebesgue exponent where one expects the behavior of eigenfunctions which maximize this quotient to change.

Thus, we do not expect '$|\mathcal{L}_x| > 0$' to be a relevant mechanism in producing large $L^p$ norms below the critical exponent. In some sense, existence of stable elliptic closed geodesics is more likely to be involved. We plan to study the problem elsewhere.

Department of Mathematics, Johns Hopkins University, Baltimore, MD 21218, USA